\documentclass[11pt]{article} 

\usepackage{pdfsync}
\usepackage[centertags]{amsmath}
\usepackage{amsfonts}
\usepackage{amssymb}
\usepackage{amsthm}
\usepackage{newlfont}
\usepackage[utf8]{inputenc}
\usepackage{mathrsfs}
\usepackage{setspace}
\usepackage{fancyhdr}
\usepackage{epsfig}
\usepackage{graphicx}
\usepackage{accents}

\usepackage{makeidx}
\makeindex

\usepackage[footnotesize,bf]{caption}

\usepackage{verbatim} 

\DeclareMathSymbol{\widehatsym}{\mathord}{largesymbols}{"62}

\newtheorem{defn}{Definition}[section]

\newtheorem{example}[defn]{Example}
\newtheorem{thm}[defn]{Theorem}
\newtheorem{prop}[defn]{Proposition}
\newtheorem{cor}[defn]{Corollary}
\newtheorem{lem}[defn]{Lemma}

\newcommand{\pf}{\noindent{\bf Proof }\mbox{   }}

\newcommand{\de}{\delta}
\newcommand{\E}{\mathcal{E}}    
\newcommand{\EV}{\mathbb{E}}    
\newcommand{\F}{\mathcal{F}}    

\newcommand{\TT}{\mathcal{T}}   
\newcommand{\Z}{\mathbb{Z}}     
\newcommand{\mc}[1]{\mathcal{#1}}

\newcommand{\ZK}{\mathbb{Z}^2_K} 
\newcommand{\bd}{\partial}	    

\newcommand{\inv}{^{-1}}

\newcommand{\ie}{\emph{i.e.}, }
\newcommand{\eg}{\emph{e.g.}, }
\newcommand{\st}{such that }

\newcommand{\lc}{\left\{}
\newcommand{\rc}{\right\}}

\newcommand{\Dnsc}{D(0,n+s)^c}
\newcommand{\Dnc}{D(0,n)^c}
\newcommand{\Dn}{D(0,n)}
\newcommand{\Dns}{D(0,n+s)}

\newcommand{\hp}{\hat{\pi}_K}
\newcommand{\hDnc}{\hp(\Dnc_K)}
\newcommand{\hDnsc}{\hp(\Dnsc_K)}
\newcommand{\hDn}{\hp(\Dn)}
\newcommand{\hDns}{\hp(\Dns)}
\newcommand{\hdDns}{\hp(\bd D(0,n)_s)}
\newcommand{\hdDnsc}{\hp((\bd D(0,n)_s)^c_K)}
\newcommand{\hDRc}{\hp(D(0,R)^c_K)}
\newcommand{\hDR}{\hp(D(0,R))}

\newcommand{\hDr}{\hp(D(0,r))}

\newcommand{\hTDnc}{T_{\hDnc}}

\newcommand{\hTDn}{T_{\hp(D(0,n))}}
\newcommand{\hTdDns}{T_{\hp(\bd D(0,n)_s)}}
\newcommand{\hTDr}{T_{\hp(D(0,r))}}

\newcommand{\hTDRc}{T_{\hp(D(0,R)^c_K)}}

\newcommand{\hGDnc}{\hat{G}_{\hp(\Dnc_K)}}

\newcommand{\hGDn}{\hat{G}_{\hp(D(0,n))}}

\newcommand{\Dr}{\hp(\bd D(x,r)_s)}
\newcommand{\Drc}{\hp((\bd D(x,r)^c_s)_K)}
\newcommand{\DRc}{\hp(D(x,R)^c_K)}
\newcommand{\DRcc}{\hp(D(x,R))}

\newcommand{\gbar}{\overline{\gamma}}

\def\midformat{
\setlength{\itemsep}{0pt} \setlength{\parindent}{0mm}
\setlength{\parskip}{0.12in} \setlength{\textheight}{180mm} 
\setlength{\textwidth}{150mm} \setlength{\evensidemargin}{0in}
\setlength{\oddsidemargin}{0in} \setlength{\topmargin}{0in}
\setlength{\hoffset}{1.0cm} \setlength{\voffset}{0.0cm}
\setlength{\headheight}{15pt} \setlength{\headsep}{.5in}
\setlength{\headwidth}{150mm} } \midformat

\newtheoremstyle{theorem}
 {}
 {}
 {\itshape}
 {}
 {\ttfamily}
 {.}
 {.5em}
 {}

\newtheoremstyle{plaintext}
 {}
 {}
 {\upshape}
 {}
 {\ttfamily}
 {.}
 {.5em}
 {}

\theoremstyle{theorem}

\theoremstyle{plaintext}

\numberwithin{equation}{section}
\begin{document}
\makeatother
\pagenumbering{roman} \thispagestyle{empty}

\title{Late Points and Cover Times\\of Projections of Planar Symmetric Random Walks \\on the Lattice Torus}
\author{Michael Carlisle\\
Baruch College, CUNY\footnote{michael.carlisle@baruch.cuny.edu}}
\date{February 20, 2013}
\maketitle

\begin{abstract}
We examine the sets of \emph{late points} 
of a symmetric random walk on $\Z^2$ projected onto the torus $\Z^2_K$, culminating in 
a limit theorem 
 for the \emph{cover time} 
of the toral random walk.
This extends the work done for the simple random walk in \cite{DPRZ2006} to a large class of random walks 
projected onto the lattice torus. The approach uses comparisons between planar and toral hitting times and distributions on annuli, and uses only random walk methods. 
\end{abstract}



\oddsidemargin 0.0in \textwidth 6.0in \textheight 8.5in

\singlespacing
\pagenumbering{arabic}
\pagestyle{fancy}
\renewcommand{\headrulewidth}{0.0pt}
\lhead{} \chead{} \rhead{\thepage}
\fancyhead[RO]{\thepage} \fancyfoot{}

\section{Introduction} \label{ch:Intro}


Wilf, in \cite{Wilf}, describes watching a simple random walk on a computer screen, where, on each time step, a dark pixel turns (and remains) bright if the walk visits it for the first time. How many steps, he wonders, will it take on average for the nearest neighbor walk's path (wrapping at the edges of the screen, making a discrete two-dimensional torus) to fill the screen? He refers to this as the ``white screen time'' problem.


He gives solutions of the white screen problem for the one dimensional path and cycle, and the complete graph $\mc{K}_n$ (known as the \emph{coupon collector}'s problem), and refers to research related to the white screen problem under the name of \emph{covering times}.  Leaving the original problem unresolved, Wilf points to a 1989 work of Zuckerman which gives bounds on the two-dimensional square lattice torus $\Z^2_K := \Z^2 / K\Z^2$. Denoting the cover time of the graph $G$ by a random walk as $\mc{T}_{cov}(G) := \sup_{x \in G} \mc{T}(x)$, where $\mc{T}(x)$ is the first hitting time of $x$, then, for the simple random walk on $\Z^2_K$, 
\[C_1 (K \log K)^2 \leq \mc{T}_{cov}(\Z^2_K) \leq C_2 (K \log K)^2\]
for some positive constants $C_1, C_2$.

Over the course of the next 20 years, closely related problems were solved by Aldous (\cite{AF}), Dembo, Peres, Rosen, \& Zeitouni (\cite{DPRZ2001}, \cite{DPRZ2004}, \cite{DPRZ2006}), Lawler (\cite{LawInt}, \cite{LawCov}, \cite{LawPol}), Rosen (\cite{RosenET}), and Rosen \& Bass (\cite{BRFreq}). This paper builds on these works to examine the structure of the so-called \emph{late points} (those not hit until ``soon'' before the cover time) which Wilf refers to as allowing the viewer of a slowly-filling white screen to ``safely go read \emph{War and Peace} without missing any action.''

We are interested in the number of late points on the square torus $\Z^2_K$ for large, increasing $K$, and will investigate this for a class of projected planar lattice, \ie $\Z^2$, random walks $S_t = S_0 + \sum_{j=0}^t X_j$, \index{s@$S_t$} for $X = \{X_j\}_{j \in \mathbb{N} \cup \{0\}}$ \index{x@$X_j$} with the following properties: $S$ is symmetric recurrent, $X_1$ has finite covariance matrix equal to a scalar times the identity, \ie $\Gamma := cov(X_1) = cI$, $c>0$, and $X$ is strongly aperiodic.\footnote{\cite{BRFreq} requires the covariance matrix of $X_1$ to be equal to $\frac{1}{2} I$, but this is a convenience for three technical points (on pages 9, 12, and 42), relating only to rotations. 
It is worthy (if not elementary) to note that the simple random walk on $\Z^d$'s $X_1$ covariance matrix is cov$(X_1) = \frac{1}{d}I$. If $K$ is odd, this walk projects to a strongly aperiodic simple random walk on $\Z^d_K$.} 
$X_1$ has, for some $\beta>0$ and $M := 4 + 2\beta$, \index{M@$M$}
\begin{equation} \label{eqn:moments}
\EV|X_1|^{M} = \sum_{x \in \Z^2} |x|^{M} p_1(x) < \infty,
\end{equation}
where, as usual in the literature,
\begin{align*}
p_1(x,y) = p_1(y-x) = P^x(X_1 = y)
\end{align*}
\index{p@$p_1$} is the one-step transition probability. 
The random walk methods used in this paper require $M > 4$; this seems to be necessary for certain Harnack inequalities which we develop (whereas, in \cite{BRFreq}, $M = 3 + 2\beta$ sufficed for frequent points on the plane).

$X$ satisfies {\bf Condition A}\footnote{Bolded terms are terms that were introduced in a paper descended from \cite{DPRZ2004} (including this author's papers), and italicized terms are well-known in the literature on random walks. 
} \index{Condition A} if either $p_1$ has bounded support, or, from any point ``just outside'' a disc, we will enter the disc with positive probability; \ie for any $s \leq n$, for large enough $n$, 
\begin{equation} \label{eqn:ConditionA}
\inf_{y: n \leq |y| < n+s} \sum_{z \in D(x,n)} p_1(y,z)
 = \inf_{y \in \bd D(x,n)_s} P^y(X_1 \in D(x,n)) \geq c e^{-\beta s^{1/4}},
\end{equation}
where the (Euclidean) $s$-annulus around the disc $D(x,n)$ (also called an \emph{$x$-band}) \index{x@$x$-band} is defined as \index{$\bd D(x,n)_s$}
\begin{equation} \label{eqn:s-band}
\bd D(x,n)_s := D(x,n+s) \setminus D(x,n).
\end{equation}
In particular, if $X_1$ has infinite range, then for any $y \in \bd D(0,n)_s$, there exists $x \in D(0,n)$ such that $p_1(y,x) > 0$.

We will switch between the planar and toral representations of the random walk and corresponding stopping times, hitting distributions, etc. 
Define the projections, \index{pi@$\pi_K$} \index{pihat@$\hp$} for $x = (x_1, x_2) \in \Z^2$, by 
\[\begin{array}{ll}
\pi_K: \Z^2 \to [-K/2,K/2)^2 \cap \Z^2, \\
\pi_K(x) = \left((x_1 + \lfloor \frac{K}{2} \rfloor) (\text{mod } K) - \lfloor \frac{K}{2} \rfloor, (x_2 + \lfloor \frac{K}{2} \rfloor) (\text{mod } K) - \lfloor \frac{K}{2} \rfloor\right); \\
 \hp: \Z^2 \to \Z^2_K, 
\,\,\,\, \hp(x) = (\pi_K x) + (K \Z)^2.
\end{array}\]
(For example, if $x = (-12,6)$ and $K = 11$, then $\pi_{11}(\Z^2) = \{-5,\ldots,5\}^2$, $\pi_{11}(x) = (-1,-5)$, and $\hat \pi_{11}(x) = (-1,-5) + (11\Z)^2$.)

\begin{figure}[!ht]
  \centering
    \includegraphics[width=6in]{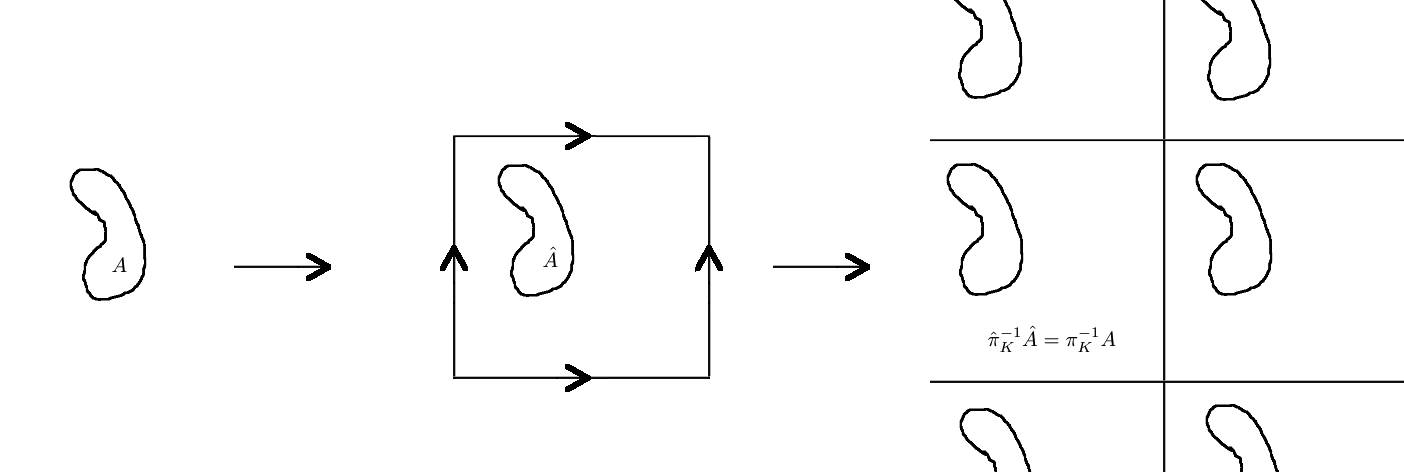}
  \parbox{4in}{
  \caption{$A \to \hat{A} \to \hp\inv \hat{A} = \pi_K \inv A$}
  \label{fig:plane_torus_pullback}}
\end{figure}

We call the set of lattice points $\pi_K(\Z^2) = [-K/2,K/2)^2 \cap \Z^2$ the {\bf primary copy} \index{primary copy} in $\Z^2$,
 and for $x \in \pi_K(\Z^2)$, $\hat{x} := \hp x$ is its corresponding element in $\Z^2_K$. Any $z \in \pi_K \inv x$, $z \neq \pi_K x$, is called a {\bf copy} \index{copy} of $x$. Likewise, for a set $A \subset \Z^2$, $\hat{A} := \hp A$ is the toral projection of $A$, and the set of all copies of $A$ is \index{piinv@$\hp\inv \hat{A}$}
\[\pi_K \inv\pi_K A = \hp\inv \hat{A} := \{z \in \Z^2: z = x + (iK, jK), \,\,  i, j \in \Z, x \in A\}.\]
Figure \ref{fig:plane_torus_pullback} displays the projection of a planar set $A$ onto the torus as $\hat{A}$, and its pullback onto $\pi_K \inv A$. (If $A \subset \pi_K \Z^2$, then of course, $A = \pi_K A$.)

For a given $\hat{x} \in \ZK$, we define $x$ to be the (planar) primary copy of that element; $x := \pi_K \hp \inv \hat x$.

While $X_j$ is the $j$th step of the planar walk and $S_j$ its position at time $j$, we use $\hat{S}_j$ \index{X@$\hat{S}_j$} to denote the position of the toral walk at time $j$. The distance between two points $x,y \in \Z^2$ will be the Euclidean distance $|x-y|$; on the torus, the distance between two points $\hat{x},\hat{y} \in \Z_K^2$ \index{xhat@$\hat{x}$} will be the minimum Euclidean distance $|\hat{x}-\hat{y}| \leq K\sqrt{2}/2$. To limit the issues regarding this distance, we will restrict any discs on $\Z^2_K$ to have radius $n < K/4$ (sometimes written as a diameter constraint: $2n < K/2$).

To bound our functions, we need a precise notion of bounding distance on the lattice torus $\Z^2_K$. As in \cite{DPRZ2006}, a function $f(x)$ is said to be $O(x)$ if $f(x)/x$ is bounded, uniformly in all implicit geometry-related quantities (such as $K$). That is, $f(x) = O(x)$ if there exists a universal constant $C$ (not depending on $K$) such that $|f(x)| \leq Cx$. Thus $x = O(x)$ but $Kx$ is \emph{not} $O(x)$. A similar convention applies to $o(x)$.

Next, we will define a few terms describing the distance of a random walk step, relative to a reference disc of radius $n$ and an $s$-annulus around the disc.
A {\bf small} jump \index{jump}
 refers to a step that is short enough to possibly (but not necessarily) stay inside a disc of radius $n$ (\ie $|X_1| < 2n$).
A {\bf baby} jump 
refers to a small jump that is too short to hop over an $s$-annulus from inside a disc (\ie $|X_1| < s$).
A {\bf medium} jump 
refers to a step that is sufficiently large to hop out of a disc and past an $s$-annulus, but with magnitude strictly less than $K$, and cannot land near a toral copy of its launching point (\ie $s < |X_1| < K-2n$).
A {\bf large} jump 
is a step which, in the toral setting, would be considered ``wrapping around'' in one step (\ie $|X_1| > K-2n$).
A {\bf targeted} jump 
is a large jump which lands directly in a copy of the disc or annulus just launched from (\ie $j(K-2n) \leq |X_1| \leq j(K+2n)/\sqrt{2}$ for some $j$). These terms will aid in dealing with differences between planar and toral hitting and escape times.\footnote{We have distinguished between three types of jumps on the torus that in the planar-only case (as in \eg \cite{BRFreq}) are referred to only as large jumps.}

As in \cite{DPRZ2001}, Section 5, set $\pi_{\Gamma} := 2\pi \sqrt{\det \Gamma}$, \index{pigamma@$\pi_{\Gamma}$} and let $\alpha \in (0,1)$. (For simple random walk, $\Gamma = \frac{1}{2}I$, so $\pi_{\Gamma} = \pi$.)
We call $\hat x$ an $\alpha,K$-\emph{late point} \index{late point} of the random walk $\hat{S}$ on $\Z^2_K$ if the first hitting time of $\hat x$, $\mc{T}_K(\hat x)$, is such that $\mc{T}_K(\hat x) \geq \frac{4\alpha}{\pi_{\Gamma}}(K \log K)^2$. 
Set $\mc{L}_K(\alpha)$ \index{LK@$\mc{L}_K(\alpha)$} to be the set of $\alpha,K$-late points in $\Z^2_K$, \ie 
\[\mc{L}_K(\alpha) := \left\{ \hat x \in \Z^2_K: \frac{\mc{T}_K(\hat x)}{(K \log K)^2} \geq \frac{4\alpha}{\pi_{\Gamma}} \right\}. \]
We prove the following, generalizing \cite[Proposition 1.1]{DPRZ2006}:
\begin{thm} \label{thm:LatePoints}
For any $0 < \alpha < 1$, 
\begin{equation} \label{eq:LatePointsAlpha}
\lim_{K \to \infty} \frac{\log |\mc{L}_K(\alpha)|}{\log K} = 2(1-\alpha) \,\, \text{ in probability.}
\end{equation}
\end{thm}
As $\alpha \to 1$, a corollary of \eqref{eq:LatePointsAlpha} is that we can generalize the cover time result of \cite[Theorem 1.1]{DPRZ2004} to our class of random walks: 
\begin{cor} \label{thm:CoverTime}
\begin{equation} \label{eq:CoverTime}
\lim_{K \to \infty} \frac{\mc{T}_{cov}(\Z^2_K)}{(K \log K)^2} = \frac{4}{\pi_{\Gamma}} \,\, \text{ in probability.}
\end{equation}
\end{cor}




The paper is structured as follows. 
In Section \ref{ch:Escape}, we state results from \cite{CarVisits} about probabilities of exiting a disc, entering a disc, and entering an annulus in the plane and torus. 
With this knowledge, in Section \ref{ch:Harnack} we build fine-tuned Harnack inequalities from general results in \cite{CarHarnack} when the landing point is a nearby annulus. 
These Harnack inequalities are applied in Section \ref{ch:Excursions} to examine excursions between consecutive concentric annuli. 
Finally, in Section \ref{ch:Late} we estimate the rarity of traveling between these annuli without ever visiting their common center point (thereby deeming the path ``late'' in visiting the center).

\section{Escape, Entry Results} \label{ch:Escape}

In this section we develop the notions of hitting time and Green's function on the plane and torus, and supply relationships between the two with respect to the timing of the random walk's escape from and entry to a disc, as well as entry to an annulus, stating results from \cite{CarVisits}.

\subsection{Disc Escape}  
\label{sec:HittingEscapeTime}
 

The \emph{hitting time} of a random walk to a set $A$ is defined as the stopping time \index{time@$T_A$} $T_{A} = \inf\{t \geq 0: S_t \in A\}$. Likewise, the \emph{escape time} of the walk from $A$ is the stopping time $T_{A^c}$. For a recurrent, strongly aperiodic, irreducible random walk on $\Z^2$, $T_{A^c} < \infty$ a.s. We denote $T_{\hat{A}}$ to be the hitting time of $\hat{A} \subset \Z_K^2$. We will examine several relationships between planar and toral hitting times.

An immediate observation on hitting times (\eg from \cite{Spitzer}) is that, the larger the set to hit, the quicker it will be hit. If $A \subset B$, then obviously $T_B \leq T_A$. It is clear, then, that $\hp\inv \hat{A}$, as an infinite number of copies of $A \subset \Z^2$, has a quicker hitting time than just one copy of $A$. In fact, we have 
\begin{equation} \label{eq:PTtimeLaw}
T_{\pi_K \inv A} = T_{\hp\inv \hat{A}} = T_{\hat{A}}.
\end{equation}

\begin{figure}[!ht]
  \centering
    \includegraphics[width=6in]{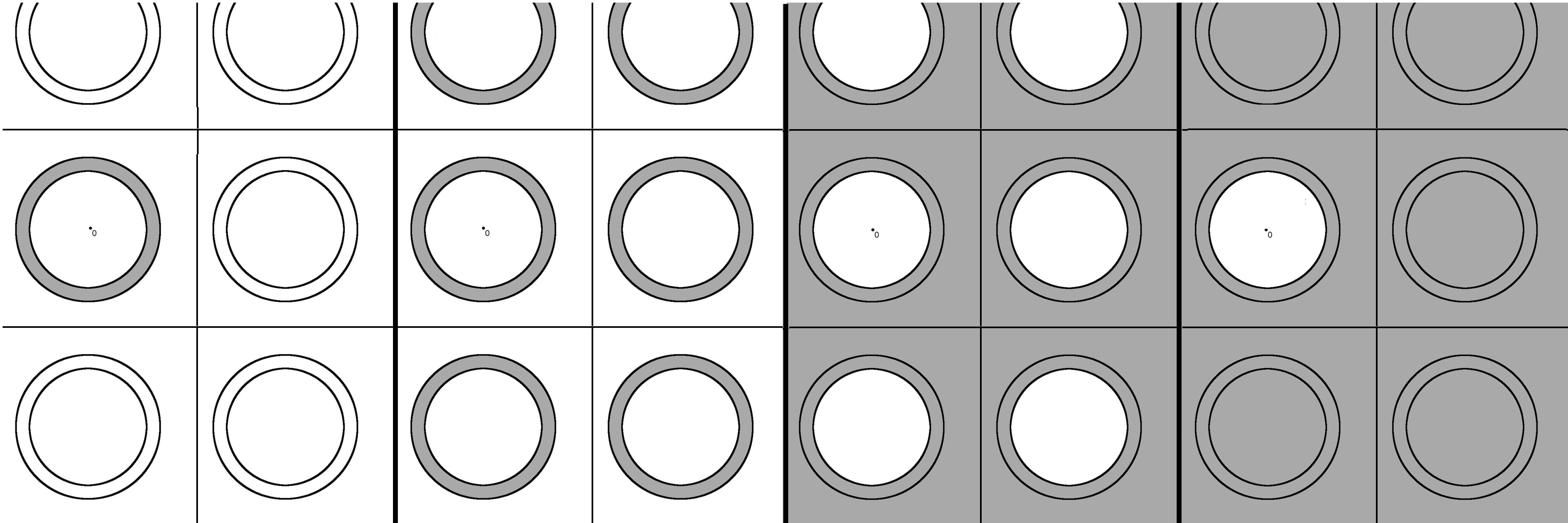}
    \parbox{6in}{
    \caption{Comparison of planar sets listed in \eqref{eq:AnnDiscSubs}, on the plane. Labeled sets are shaded.}
    \label{fig:discs_annuli}}
\end{figure}

Let $n, s$ be such that $n+s < K/4$, and $D(0,n) = \pi_K \Dn$ the primary copy of $\Dn \subset \Z^2$. Define the primary copy's portion of the complement of $\Dn$ to be $\Dnc_K := \Dnc \cap \pi_K \Z^2$. \index{DncK@$\Dnc_K$}
\eqref{eq:AnnDiscSubs} and Figure \ref{fig:discs_annuli} describe the nestedness of sets from the planar annulus $\bd D(0,n)_s$ up to the planar disc complement $\Dnc$: 
\begin{align}
\bd D(0,n)_s & \subset \pi_K \inv (\bd D(0,n)_s) = \hp\inv \hdDns \notag\\
 & \subset \hp\inv \hDnc = \pi_K \inv (\Dnc_K) \subset D(0,n)^c. \label{eq:AnnDiscSubs}
\end{align}

By \eqref{eq:PTtimeLaw}, \eqref{eq:AnnDiscSubs} yields, starting at any $x \in D(0,n)$, the \emph{disc escape time inequalities}
\begin{align} 
T_{\bd D(0,n)_s} & \geq T_{\pi_K \inv \bd D(0,n)_s} = T_{\hp\inv \hdDns} \notag \\
 & \geq T_{\hp\inv \hDnc} = T_{\pi_K \inv (\Dnc_K)} \geq T_{D(0,n)^c} \geq 1. \label{eq:HittingTimeComp1}
\end{align}

We shall take planar starting points from the primary copy ($x = \pi_K x$). The probabilities of these inequalities being strict (\eg $P^x(T_{D(0,n)^c} < \hTDnc)$) and the means of the stopping times will be of interest to us. We start with estimating the mean of the planar escape time from $D(0,n)$ (which improves on \cite[Prop. 6.2.6]{LawSRW}), and then use this probability to estimate the toral escape time from $\hDn$.

\begin{lem} \label{lem:EscapeDisc} 
Let $S_t = S_0 + \sum_{j=1}^t X_j$ be a random walk in $\Z^2$ with $E|X_1|^2 < \infty$, and covariance matrix $\Gamma$ such that $tr(\Gamma) = \gamma^2 > 0$. Then,
uniformly for $x \in D(0,n)$, and for sufficiently large $n$, \index{expected@$\EV^x(T_{D(0,n)^c})$}
\begin{equation} \label{eqn:EscapeDiscExp}
\frac{n^2 - |x|^2}{\gamma^2} \leq \EV^x(T_{D(0,n)^c}) \leq \frac{n^2 - |x|^2}{\gamma^2} + 2n + 1.
\end{equation}
\end{lem} 

\pf See \cite[Lemma 2.1]{CarVisits}. \qed

For $\Gamma = cI$, $\gamma^2 = 2c$ and so \eqref{eqn:EscapeDiscExp} becomes\footnote{For simple random walk on $\Z^2$, $c = 1/2$, which yields \cite[(2.3)]{DPRZ2006}.} 
\begin{equation} \label{eqn:EscapeDiscExpCalc}
\frac{n^2 - |x|^2}{2c} \leq \EV^x(T_{D(0,n)^c}) \leq \frac{n^2 - |x|^2}{2c} + 2n + 1.
\end{equation}

We define the \emph{Green's function} \index{Green@$G_{t^*}$} for two points $x, y$, as the expected number of visits to $y$, starting from $x$, up to the fixed time $t^*$:
\begin{equation} \label{eq:GreenDefn}
G_{t^*}(x,y) := \EV^x\bigg[\sum_{j=0}^{t^*} 1_{\{S_j=y\}}\bigg] = \sum_{j=0}^{\infty} P^x(S_j = y; j < t^*).
\end{equation}
Spitzer, in \cite{Spitzer}, similarly defines the \emph{truncated Green's function}, \index{Green@$G_A$} for $x, y \in A$ of a random walk from $x$ to $y$ before exiting $A$ as the total expected number of visits to $y$, starting from $x$:
\begin{equation} \label{eq:GreenDefn}
G_A(x,y) := \EV^x\bigg[\sum_{j=0}^{\infty} 1_{\{S_j=y; j<T_{A^c}\}}\bigg] = \sum_{j=0}^{\infty} P^x(S_j = y; j < T_{A^c})
\end{equation}
and 0 if $x$ or $y \not \in A$. (Since the walk is recurrent and aperiodic, there is no ``all-time'' Green's function to count the total number of visits to $x$ from $j=0$ to $\infty$.) An elementary result for any random walk (found, for example, in \cite{Spitzer}, or \cite[Sect. 1.5]{LawInt}) is that, for $x, y \in A \subset B$, there are more possible visits inside $B$ than inside $A$:
\begin{equation} \label{eq:GreenComp}
 G_A(x,y) \leq G_B(x,y).
\end{equation}
Also of interest is the expected hitting time identity 
\begin{equation} \label{eq:GreenHit}
\EV^x(T_{A^c}) = \sum_{z \in A} G_A(x,z).
\end{equation}
Starting at a point $x \in A^c$, the \emph{hitting distribution} \index{hitting@$H_A$} of $A$ 
 is defined as
\[H_A(x,y) := P^x(S_{T_A} = y).\]
The \emph{last exit decomposition} \index{last exit decomposition} of a hitting distribution is based on the Green's function: for $A$ a proper subset of $\Z^2$, $x \in A^c$, $y \in A$, 
\begin{equation} \label{eq:LastExit}
H_{A}(x,y) = \sum_{z \in A^c} G_{A^c}(x,z) p_1(z,y).
\end{equation}
An immediate result follows from \eqref{eq:GreenComp}:  
If $y \in A \subset B$, then for $x \in B^c \subset A^c$, we have by \eqref{eq:GreenComp} the monotonicity result
\begin{equation} \label{eq:HittingComp}
\begin{array}{lll}
H_{A}(x,y) & = & \sum_{z \in A^c} G_{A^c}(x,z) p_1(z,y)\\
 & \geq & \sum_{z \in B^c} G_{B^c}(x,z) p_1(z,y) = H_{B}(x,y)
\end{array}
\end{equation}
and the subset hitting time relations (assuming a recurrent random walk) 
\begin{align} 
P^x(T_A = T_B) & = \sum_{z \in A} H_B(x,z); \notag\\
P^x(T_A \neq T_B) & = P^x(T_A > T_B) = \sum_{z \in B \setminus A} H_B(x,z) \label{eq:SubsetHittingComp}
\end{align}
which we will revisit in Section \ref{ch:Annulus}.

By Markov's inequality, large jumps are rare: if $C_M = \EV(|X_1|^M) < \infty$, then since $2n < K/2$, 
\begin{equation} \label{eqn:LargeJumpEstimate}
P(|X_1|>K-2n) \leq \frac{C_M}{(K-2n)^M} < \frac{2^M C_M}{K^M} = O(K^{-M}).
\end{equation}
Recall that, when given a toral element $\hat{x} \in \ZK$, we define $x$ to be the (planar) primary copy of that element; $x := \pi_K \hp \inv \hat x$. A toral step $\hat{x} \to \hat{y}$ must take into account large jumps that, on the plane, would land on a copy of $y$ (\ie in $\hp\inv \hat{y}$). All of these positions, together, are a small addition to the planar jump probability. By \eqref{eqn:LargeJumpEstimate} we have, for $\hat{x}, \hat{y} \in \Z^2_K$, the targeted jump estimate \index{phat@$\hat{p}_1$}
\begin{align} \label{eqn:TargetedJumpEstimate}
\hat{p}_1(\hat x, \hat y) := P^{\hat{x}}(\hat S_1 = \hat{y}) & = P^{x}(S_1 = y) + P^{x}\left(|X_1|>K-2n ; S_1 \in \hp\inv\hat{y} \setminus \{y\}\right) \notag\\
 & \leq p_1(x, y) + O(K^{-M}).
\end{align}

By \eqref{eq:LastExit}, \eqref{eqn:LargeJumpEstimate}, and then \eqref{eqn:EscapeDiscExp} and \eqref{eq:GreenHit}, for some $c < \infty$ and any $x\in D(0,n)$, 
\begin{align}
P^x(\hTDnc>T_{D(0,n)^c}) & = \sum_{z \in \left(\hp\inv\hDn \,\setminus\, D(0,n)\right)} \sum_{y \in D(0,n)} G_{D(0,n)}(x,y) p_1(y,z) \notag\\
 & \leq c K^{-M} \sum_{y \in D(0,n)} G_{D(0,n)}(x,y) = O(K^{-M} n^2). \label{eqn:ProbTorusExit}
\end{align}

We now find that the mean of the disc escape time on the torus is larger than on the plane, but only by a small factor (induced by the rarity of targeted jumps). \index{expectedT@$\EV^{\hat x}(\hTDnc)$}

\begin{lem} \label{lem:EscapeBounds} 
For $n < K/4$, $x \in D(0,n)$, and $n$ and $K$ sufficiently large,  
\begin{equation} \label{eqn:Escape}
\EV^{\hat x}[\hTDnc] \leq \EV^x[T_{D(0,n)^c}] + O(K^{-M} n^2) \max_{y\in D(0,n)}\EV^y[T_{D(0,n)^c}].
\end{equation}
\end{lem}

\pf See \cite[Lemma 2.2]{CarVisits}. \qed

\begin{example} \label{eg:WorstCaseTiming}
Let $A=D(0,\sqrt{2}) = \{0, +e_1, -e_1, +e_2, -e_2\} \subset \Z^2$, where $e_i$ is the $i$th unit vector in $\Z^2$, and $K$ odd and fixed. Let $X$ be the symmetric random walk on $\Z^2$ starting at $X_0=0$ defined by the probabilities
\[p_1(K^j e_i) = P^0(X_1= K^j e_i) = \frac{1}{4}e^{-\lambda}\frac{\lambda^j}{j!}, j=0,1,2,\ldots; \,\,\, i=1,2.\]
$\frac{\log |X_1|}{\log K}$ is a Poisson random variable with parameter $\lambda$, and moving any of the four primary lattice directions is equally likely. $S_t$ is strongly aperiodic recurrent and has infinite range, $E(|X_1|^m) < \infty$ for all $m < \infty$ (and, in particular, cov$(|X_1|) = \Gamma = \frac{1}{2} e^{(K^2 - 1)\lambda} I$), and every large jump causes a landing in a new copy of $A$. The only way to escape $\pi_K \inv A = \hp \inv \hat{A}$ is a step of size $K^0=1$.
\end{example}

Computational bounds on $\EV^{\hat x}(\hTDnc)$, by \eqref{eqn:Escape} and \eqref{eqn:EscapeDiscExp}, are 
\begin{equation} \label{eqn:EscapeDiscExpTorus}
\frac{n^2 - |x|^2}{\gamma^2} \leq \EV^{\hat x}(\hTDnc) \leq \frac{n^2 - |x|^2}{\gamma^2} + 2n + 1 + O(K^{-M} n^4).
\end{equation}

\begin{example} \label{eg:LazySRW}
Define the $\varepsilon$-lazy simple random walk on $\Z^d$, for $0 \leq \varepsilon < 1$, to be the walk with steps 
$p_1(e_j) = p_1(-e_j) = \frac{1 - \varepsilon}{2d}$, $j=1,...,d$; $p_1(0) = \varepsilon$,
\ie the walk stands still for a step with probability $\varepsilon$, and acts ``simply'' otherwise. Then $\Gamma = \left(\frac{1-\varepsilon}{d}\right)I$, and so for $d=2$, $\EV^{\hat x}(\hTDnc) = \frac{n^2 - |x|^2}{1-\varepsilon} + O(n)$.
\end{example}

We will next see that, from inside a disc, the probability of hitting zero before escaping is nearly the same on the torus as on the plane. Recall that, for $\hat{x} \in \ZK$, $x := \pi_K \hp \inv \hat x$.

\begin{lem} \label{lem:HitZeroFirst}
For all $\hat{x} \in \hDn$ and $n$ sufficiently large with $2n < K/2$,
\begin{equation} \label{eq:HitZeroFirst}
P^{\hat{x}}(T_{\hat{0}} < \hTDnc) = P^{x}(T_0 < T_{D(0,n)^c}) + O(K^{-M} n^2).
\end{equation}
\end{lem}

\pf See \cite[Lemma 2.3]{CarVisits}. \qed

Finally, we calculate bounds for hitting time probabilities of a small disc around zero before escaping the $n$-disc. Let $\rho(\hat{x}) := n-|\hat{x}|$ be the distance between $\hat{x}$ and $\hp(D(0,n))$. 

\begin{lem}
Let $0 < \delta < \varepsilon < 1$. Then there exist $0 < c_1 < c_2 < \infty$ such that for all $\hat x \in \hDn \setminus \hp(D(0,\varepsilon n))$, for $n$ sufficiently large,
\begin{equation} \label{eq:InnerHitBounds}
c_1 \frac{\rho(\hat x) \vee 1}{n} \leq P^{\hat x}(T_{\hp(D(0,\delta n))} < \hTDnc) \leq c_2 \frac{\rho(\hat x) \vee 1}{n}.
\end{equation}
\end{lem}

\pf See \cite[Lemma 2.4]{CarVisits}. \qed


Here we will examine \emph{internal} Green's functions on the plane (\ie from inside a disc; Green's functions external to a disc will be analyzed in Section \ref{ch:Entry}). We extend some results of \cite{LawSRW} for symmetric random walks on $\Z^2$ to projections of these random walks onto $\Z_K^2$.

We define the Green's function in the usual way for $\hat x, \hat y \in \hp(A) = \hat{A} \in \ZK$ to be, in comparison to \eqref{eq:GreenDefn},  
\begin{equation} \label{eqn:Green}
\hat{G}_{\hp(A)}(\hat x,\hat y) := \sum_{j=0}^{\infty} P^{\hat x}(\hat S_j = \hat y; j < T_{\hp(A^c_K)}) 
\end{equation}
and 0 else. In the planar case, the stopping time $T_{A^c}$ for a bounded set $A$ has a clear meaning, as a sufficiently large jump (one with magnitude greater than the diameter of $A$, for example) will certainly exit $A$. Jumps targeting $A$ land, in $\Z^2$, in $\pi_K\inv A = \hp\inv\hat{A}$; on $\ZK$, they land in $\hat{A}$. This means that planar estimates 
must be adjusted to reach similar results on the torally-projected walk. 

Please note that \eqref{eqn:Green} is different from the planar Green's function on the periodic planar set $\pi_K \inv A$: 
\begin{equation} \label{eqn:GreenPT}
G_{\pi_K \inv A}(x,y) := \sum_{j=0}^{\infty} P^x(S_j = y; j < T_{\pi_K \inv (A^c_K)}), \,\, x, y \in \pi_K \inv A.
\end{equation}
We will explore this distinction in Section \ref{ch:Entry}.

Note that $S_j \in \hp\inv\hat{S}_j$ for every $j$. By \eqref{eq:HittingTimeComp1} it is clear that planar escape happens at or before toral escape. Hence, the number of planar visits is less than or equal to the number of toral visits; for any $x,y \in A \subset \pi_K \Z^2$, 
\begin{align} 
G_A(x,y) & = \sum_{j=0}^{\infty} P^x(S_j = y; \, j < T_{A^c}) \notag\\
 & = \sum_{j=0}^{\infty} P^x(S_j \in \pi_K\inv y; \, j < T_{A^c}) 
 = \sum_{j=0}^{\infty} P^{\hat x}(\hat{S}_j = \hat{y}; \, j < T_{A^c}) \label{eqn:IntGreenIneq}\\
 & \leq \sum_{j=0}^{\infty} P^{\hat x}(\hat{S}_j = \hat{y}; \, j < T_{\hp(A^c_K)}) = \hat{G}_{\hp(A)}(\hat x, \hat y), \notag
\end{align}
where equality occurs between the first and second lines because, of all the copies of $y$ in $\pi_K\inv y$, only the primary copy $y = \pi_K y$ can be hit before the planar escape time $T_{A^c}$.

We start by giving bounds on the number of visits to $\hat{0}$ before escaping a disc. 
\begin{lem} \label{lem:GreenZero}
For $n$ sufficiently large (with $2n < K/2$),
\begin{equation} \label{eqn:GreenZero}
\hGDn(\hat 0,\hat 0) = G_{D(0,n)}(0,0)[1 + O(K^{-M} n^2)].
\end{equation}
\end{lem}

\pf See \cite[Lemma 2.5]{CarVisits}. \qed

Define the \emph{potential kernel} for $X$ on $\Z^2$ as follows: for $x \in \Z^2$, \index{axy@$a(x)$} 
\begin{align}
a(x) := \lim_{n \to \infty} \sum_{j=0}^n [p_j(0) - p_j(x)]. \label{eq:axyDefn}
\end{align}
Combining the generality of rotation of \cite[Ch. III, Sec. 12, P3]{Spitzer} and \cite[Theorem 4.4.6]{LawSRW} and the infinite-range argument of \cite[Prop. 9.2]{BRFreq} gives, for covariance matrix $\Gamma$ and norm $\mc{J}^*(x) := |x \cdot \Gamma^{-1} x|$, as $|x| \to \infty$, 
\begin{align}
a(x) = \frac{2}{\pi_{\Gamma}} \log\mc{J}^*(x) + C(p_1) + o(|x|^{-1}), \label{eq:axyGeneralCalc}
\end{align}

where $C(p_1)$ is a constant depending on $p_1$ but not $x$, and $\pi_{\Gamma} = 2\pi\sqrt{\det \Gamma}$. For $\Gamma = cI$, this reduces to 
\begin{align}
a(x) & = \frac{1}{c \pi} \log\left(\frac{|x|}{\sqrt{c}}\right) + C(p_1) + o(|x|^{-1}) \notag\\
 & = \frac{1}{c \pi} \log|x| + C'(p_1) + o(|x|^{-1}) , \label{eq:axyCI}
\end{align}
where $C'(p_1) = C(p_1) - \frac{1}{2c\pi}\log c$. For simple random walk on $\Z^2$, $c=\frac{1}{2}$, and so this is, from \cite[Theorem 4.4.4]{LawSRW}, 
\begin{align}
a(x) & = \frac{2}{\pi} \log|x| + \frac{2\gamma + \log 8}{\pi} + o(|x|^{-1}), \label{eq:axy}
\end{align}
where $\gamma$ is Euler's constant. From here on, we will write \eqref{eq:axyCI} with the form 
\begin{align}
a(x) = \frac{2}{\pi_{\Gamma}} \log|x| + C'(p_1) + o(|x|^{-1}). \label{eq:axyCIgen}
\end{align}
By the argument in \cite[(2.8)-(2.12)]{BRFreq} (which calculates the overshoot estimate of $O(n^{-1/4})$ mentioned in the note after \cite[Prop. 6.3.1]{LawSRW}), and using \eqref{eq:axyCIgen}, we get a computational result for \eqref{eqn:GreenZero} if $\Gamma = cI$: 
\begin{align}
G_{D(0,n)}(0,0) & = \frac{2}{\pi_{\Gamma}}\log n + C' + O(n^{-1/4}) \label{eq:FP(2.13)} 
\end{align}
which implies the toral Green's function 
\begin{align}
\implies \hGDn(\hat 0,\hat 0) & = G_{D(0,n)}(0,0)(1 + O(K^{-M} n^2)) \notag\\
 & = \left(\frac{2}{\pi_{\Gamma}}\log n + C' + O(n^{-1/4})\right)(1 + O(K^{-M} n^2)) \notag\\
 & = \frac{2}{\pi_{\Gamma}}\log n + C' + O(n^{-1/4}). \label{eq:GreenZeroTorusVal}
\end{align}
For $x, y \in \Z^2$ such that $|x| \ll |y|$, we have, by a Taylor expansion around $y$, 
\begin{align} 
\log|y-x| & = \log|y| + O\left( \frac{|x|}{|y|} \right). \label{eq:logEstimate}
\end{align}
In particular, if $x \in D(0,2r)$ and $y \in D(0,R/2)^c$, with $R = 4mr$, we have 
\begin{align} 
\log|y-x| & = \log|y| + O\left( m^{-1} \right). \label{eq:logEstimate2}
\end{align}
Note that \eqref{eq:logEstimate} and \eqref{eq:logEstimate2} hold in the toral case without adjustment.

Let $\eta = \inf\{ t \geq 1: S_t \in \{0\} \cup D(0,n)^c\}$. Then, following the argument of \cite[(2.14)-(2.15)]{BRFreq}, since $a(x)$ is harmonic with respect to $p$, $a(S_{t \land \eta})$ is a bounded martingale. Hence, $|a(S_{t \land \eta})|^2$ is a submartingale, so $\EV|a(S_{t \land \eta})|^2 \leq \EV|a(S_{\eta})|^2 < \infty$, meaning $\{a(S_{t \land \eta})\}$ are uniformly integrable. Hence, by the optional stopping and bounded convergence theorems, \eqref{eq:axyCIgen}, and \eqref{eq:logEstimate2}, 
\begin{align*}
a(x) & = \lim_{t \to \infty} \EV^x(a(S_{t \land \eta})) = \EV^x(a(S_{\eta})) = \EV^x(a(S_{\eta}); \, S_{\eta} \neq 0) \\
 & = \sum_{y \in \bd D(0,n)_{n^{3/4}}} a(y) P^x(S_{\eta} = y) + \sum_{y \in D(0,n+n^{3/4})^c} a(y) P^x(S_{\eta} = y) \\
 & = \left( \frac{2}{\pi_{\Gamma}} \log n + C'(p_1) + o(|x|^{-1}) + O(n^{-1/4})\right) P^x(S_{\eta} \neq 0) + O(n^{-1/4}), 
\end{align*}
which, combining the error terms into $O(|x|^{-1/4})$, matches \cite[Prop. 6.4.3]{LawSRW}: 
\begin{align}
P^x&(T_0 < T_{D(0,n)^c}) = P^x(S_{\eta} = 0) = 1 - \frac{a(x) - O(n^{-1/4})}{\frac{2}{\pi_{\Gamma}} \log n + C' + O(|x|)^{-1/4}} \label{eq:FP(2.16)}\\
 & = 1 - \frac{\frac{2}{\pi_{\Gamma}} \log |x| + C' + O(|x|^{-1/4})}{\frac{2}{\pi_{\Gamma}} \log n + C' + O(n^{-1/4})} 
 = \left( \frac{\log (n/|x|) + O(|x|^{-1/4})}{\log n}\right)(1 + O((\log n)^{-1})). \notag
\end{align}
With \eqref{eq:HitZeroFirst}, 
we move this to the torus: 
\begin{align}
P^{\hat x}(T_{\hat{0}} < \hTDnc) & = \frac{\log(n/|\hat{x}|) + O(|\hat{x}|^{-1/4})}{\log(n)}\bigg(1 + O((\log n)^{-1})\bigg)   + O(K^{-M} n^2) \notag\\
 & = \frac{\log(n/|\hat{x}|) + O(|\hat{x}|^{-1/4})}{\log(n)}\bigg(1 + O((\log n)^{-1})\bigg). \label{eq:ProbZeroBeforeDisc}
\end{align}

Next, we examine $\hat{x} \in \hDR \setminus \hDr$. By the fact that a large targeted jump may land a planar walk into $\hp\inv \hDr \setminus D(0,r)$ (the set of any copy of $D(0,r)$ that is not the primary copy), we may transfer the planar results \cite[(2.20), (2.21)]{BRFreq} 
\begin{align} 
P^x(T_{D(0,r)} > T_{D(0,R)^c}) & = \frac{\log(|x|/r) + O(r^{-1/4})}{\log(R/r)} \label{eq:FP(2.20)} \\
P^x(T_{D(0,r)} < T_{D(0,R)^c}) & = \frac{\log(R/|x|) + O(r^{-1/4})}{\log(R/r)} \label{eq:FP(2.21)}
\end{align}
uniformly for $r < |x| < R$ to the toral results 
\begin{align}
P^{\hat x}(\hTDr > \hTDRc) & = \frac{\log(|\hat{x}|/r) + O(r^{-1/4})}{\log(R/r)} + O(K^{-M} R^2) \notag\\
 & = \frac{\log(|\hat{x}|/r) + O(r^{-1/4})}{\log(R/r)} \label{eq:ToralGamblersSuccess}  \\
P^{\hat x}(\hTDr < \hTDRc) & = \frac{\log(R/|\hat{x}|) + O(r^{-1/4})}{\log(R/r)} + O(K^{-M} R^2) \notag\\
 & = \frac{\log(R/|\hat{x}|) + O(r^{-1/4})}{\log(R/r)}. \label{eq:ToralGamblersRuin}
\end{align}

The strong Markov property applied at $T_0$ gives us the planar equality
\begin{equation} \label{eq:GreenXZeroPlanar}
G_{D(0,n)}(x,0) = P^x(T_0 < T_{D(0,n)^c}) \, G_{D(0,n)}(0,0), 
\end{equation}
which implies $G_{D(0,n)}(x,0) \leq G_{D(0,n)}(0,0)$ for any $x \in D(0,n)$. This equality has a clear analog on the torus, by applying the strong Markov property at $T_{\hat 0}$, for any $\hat x \in \hDn$,
\begin{equation} \label{eq:GreenTorusXZero}
\hGDn(\hat x,\hat 0) = P^{\hat{x}}(T_{\hat{0}} < \hTDnc) \, \hGDn(\hat 0,\hat 0).
\end{equation}
By \eqref{eq:GreenXZeroPlanar}, \eqref{eq:FP(2.13)}, \eqref{eq:GreenZeroTorusVal}, \eqref{eq:GreenTorusXZero}, \eqref{eq:FP(2.16)}, and \eqref{eq:ProbZeroBeforeDisc}, we get as corollaries calculations and bounds for $G_{D(0,n)}(x,0)$,  $\hGDn(\hat x,\hat 0)$, $G_{D(0,n)}(x,z)$, and $\hGDn(\hat x,\hat z)$: for $x \in \Dn$ and $\hat x \in \hDn$, 
for some $C = C(p) < \infty$, 
\begin{align} 
G_{D(0,n)}(x,0) & = P^x(T_0 < T_{\Dnc}) \, G_{\Dn}(0,0) \notag\\
 & = \frac{\log(n/|x|) + O(|x|^{-1/4})}{\log(n)} \left( 1 + O((\log n)^{-1} \right) \left( \frac{2}{\pi_{\Gamma}} \log n + C' + O(n^{-1/4})\right) \notag\\ 
 & = \frac{2}{\pi_{\Gamma}}\log \bigg(\frac{n}{|x|}\bigg) + C + O(|x|^{-1/4}), \label{eq:FP(2.18)} 
\end{align}
\begin{align}
\hGDn(\hat x,\hat 0) & = \frac{2}{\pi_{\Gamma}}\log \bigg(\frac{n}{|\hat x|}\bigg) + C + O(|\hat x|^{-1/4}) \label{eq:GreenXZeroTorusCalc} \\
G_{D(0,n)}(x,z) & \leq G_{D(x,2n)}(0,z-x) \leq c \log n. \label{eq:FP(2.19)} \\
\hGDn(\hat x,\hat z) & = G_{D(0,n)}(x,z) + O(K^{-M} n^2 \log n) \leq c \log n.\label{eqn:GreenXZBounds}
\end{align}

Finally, we have the following result paralleling \eqref{eq:InnerHitBounds}. Recall that $\rho(\hat{x}) = n - |\hat{x}|$.
\begin{lem} \label{lem:GreenFP2.2}
For any $0<\delta<\varepsilon<1$ we can find $0<c_{1}<c_{2}<\infty$, such that for  all $\hat x \in \hDn \setminus \hp(D(0,\varepsilon n))$, $\hat y \in \hp(D(0, \delta n))$ and all $n$ sufficiently large such that $2n < K/2$,
\begin{equation} \label{eq:GreenFP2.2G}
c_{ 1}\frac{\rho( \hat x)\lor 1}{n}\leq \hat{G}_{\hDn}(\hat y,\hat x)\leq c_{2}\frac{\rho( \hat x)\lor 1}{n}.
\end{equation}
\end{lem}

\pf See \cite[Lemma 2.6]{CarVisits}. \qed


\subsection{Disc Entry} \label{ch:Entry}

Here we will examine paths starting outside a disc. Since, on $\Z^2$, 
\begin{equation} \label{eq:AnnDiscSubs2}
\bd D(0,n)_s \subset \left\{\begin{array}{c}
\pi_K \inv \bd D(0,n)_s \\
D(0,n+s)
\end{array}\right\} \subset \pi_K \inv D(0,n+s),
\end{equation}
then starting at any $y \in \pi_K \inv (D(0,n+s)^c \cap \pi_K \Z^2)$ (as seen in Figure \ref{fig:discs_annuli}) yields the \emph{disc entrance time inequalities}
\begin{equation} \label{eq:HittingTimeComp2}
T_{\pi_K \inv D(0,n+s)} \leq \left\{\begin{array}{c}
T_{\pi_K \inv \bd D(0,n)_s} \\
T_{D(0,n+s)}
\end{array}\right\} \leq T_{\bd D(0,n)_s}.
\end{equation}

These relationships will be exploited in this and the next section.


To supplement the internal Green's functions of Section \ref{ch:Escape} are external Green's functions: those counting the number of visits to a point outside of a set before entering that set. Wlog $x$ and $D(0,n)$ are in the primary copy. We will find bounds on three different external Green's functions: 
\begin{center}
\begin{tabular}{c|c|c|c|l}
Green's function & scope & starting at & counts visits to & before...\\
\hline
$G_{D(0,n)^c}(x,y)$ & planar & $x$ & $y$ & $T_{D(0,n)}$ \\
$G_{\pi_K \inv (D(0,n)^c_K)}(x,y)$ & planar & $x$ & $y$ & $T_{\pi_K \inv D(0,n)} = \hTDn$ \\
$\hGDnc(\hat{x},\hat{y})$ & toral & $\hat{x}$ & $\hat{y}$ & $T_{\pi_K \inv D(0,n)} = \hTDn$
\end{tabular}
\end{center}

Note that, similar to \eqref{eq:GreenXZeroPlanar}, for any $x, y \in D(0,n)^c$, by the symmetry of $G_A$ and the strong Markov property at $T_x$, 
\begin{equation} \label{eq:ExternalPlanarEasyBound}
G_{D(0,n)^c}(x,y) = P^y(T_x < T_{D(0,n)}) \, G_{D(0,n)^c}(x,x),
\end{equation}
so, assuming $|x| < |y|$, we only need $G_{D(0,n)^c}(x,x)$ for an upper bound. Fix $j > 2$. By the arguments from \cite[Section 3.1]{CarVisits}, we have the following bounds for any $x,y \in \pi_K(\Dnc_K)$ \st $|x| \leq |y|$: 
\begin{eqnarray} 
G_{D(0,n)^c}(x,y) & \leq & \frac{2j}{j-2} \left[ \frac{2}{\pi_{\Gamma}} \log(2|x|) + C + O(|x|^{-1/4})\right] \leq c_j \log |x|, \label{eqn:ExtGreenIneq}\\
\hGDnc(\hat x,\hat y) & \leq & \frac{2j}{j-2} \left[ \frac{2}{\pi_{\Gamma}} \log(2|x|) + C + O(|x|^{-1/4})\right] \leq \hat{c}_j \log |\hat x|, \label{eqn:ExtGreenIneqToral}
\end{eqnarray}
where $c_j, \hat{c}_j$ depend on $j>2$, $c_j \geq \hat{c}_j$, and in the toral case, such that $|\hat x| < (\frac{K}{2})^{1/j}$ (there is no such restriction on the planar case). 
\begin{lem} \label{lem:ExternalGreenToralComplete}
For $\hat x \in \hDnc$, 
\begin{equation} \label{eqn:ExternalGreenToralComplete}
\hGDnc(\hat x,\hat x) \leq \left\{ \begin{array}{ll}
C \log |\hat x| & n < |\hat x| < \left( \frac{K}{2} \right)^{1/3} \\
C\log^{2} |\hat x| & \left( \frac{K}{2} \right)^{1/3} \leq |\hat x|.
\end{array}\right.
\end{equation}
\end{lem}

\pf See \cite[Lemma 3.1]{CarVisits}. \qed


We will now approach disc entrance times.
Our first planar result mirrors (\ref{eqn:EscapeDiscExp}), with a very different end result, which is hinted by the first passage time result for SRW on $\Z$ (in, for example, \cite{ShreveVolI}).

\begin{lem} \label{lem:EnterDisc}
For any $y \in D(0,n)^c$, \index{expected@$\EV^y(T_{D(0,n)})$}
\begin{equation} \label{eqn:EnterDiscExp}
\EV^y(T_{D(0,n)}) = \infty.
\end{equation}
\end{lem}

\pf See \cite[Lemma 3.2]{CarVisits}. \qed

Next, we find finite bounds on the expected time to enter a toral disc. 
\begin{lem}
For any $n < \frac{K}{6}$ and $\hat y \in \hDnc$, there exists $c < \infty$ such that 
\begin{equation} 
\EV^{\hat y}(\hTDn) \leq \left\{ \begin{array}{ll}
c K^2 \log \left(n\right) & n < |\hat y| < n^2 \\
c K^2 \log \left(\frac{|\hat y|}{n}\right) & n^2 \leq |\hat y| < \left( \frac{K}{2} \right)^{1/3} \label{eq:ToralDiscEntryExpectedTime}\\
c K^2 (\log |\hat y|)^2 & \left( \frac{K}{2} \right)^{1/3} \leq |\hat y|.
\end{array}\right. 
\end{equation}
Also, we have the lower bound 
\begin{equation} 
\EV^{\hat y}(\hTDn) \geq \left\{ \begin{array}{ll}
\frac{(|\hat{y}|-n)^2}{\gamma^2} & |\hat{y}| < \frac{K}{3} \\
 \frac{c(K-n)^2}{\gamma^2} & |\hat y| \geq \frac{K}{3}
\end{array} \label{eq:ToralDiscEntryExpectedTimeLower} \right. 
\end{equation}
where $\gamma^2$ is as in the proof of Lemma \ref{lem:EnterDisc}.
\end{lem}

\pf See \cite[Lemma 3.3]{CarVisits}. \qed

\eqref{eq:ToralDiscEntryExpectedTime} hints at the late points and cover time results of Section \ref{ch:Late}. We will improve on these bounds in our discussion on excursions.

\subsection{Annulus Entry} \label{ch:Annulus}

In this section we will state results from \cite{CarThree} for general Green's functions, hitting times, and hitting distributions by a symmetric recurrent random walk $X$ on a set partitioned into three pieces. We then apply these results to the partition of disc, annulus, and ``outside'' to relate our results from Sections \ref{sec:HittingEscapeTime} and \ref{ch:Entry} to the annulus. We conclude by finding tailored gambler's ruin-based probabilities and hitting distribution bounds for annuli.

\subsubsection{Bounds on a three-partitioned set} 

Let $A \sqcup B \sqcup C$ partition our sample space. We find estimates for the Green's function $G_{A \cup B}$ and the hitting time $E^x(T_C)$ for $x, y \in A \cup B$, with interest in the case where $C$ ``separates'' $A$ and $B$ in a sense (\ie the probability of jumping from $A$ to $B$, or vice versa, without hitting $C$, is small). This gives a notion for how probabilistically ``separate'' they are.

Simple lower bounds for the Green's function $G_{A \cup B}$ are obvious; to find upper bounds for these cases, we analyze excursions between $A$ and $B$ before hitting $C$.
\begin{lem} \label{lem:GreenABUB}
For $a,a' \in A$ and $b,b' \in B$, with $\theta_t$ the usual shift operators, 
\begin{align}
T^*_B := & \inf\{t > T_A: X_t \in B\} = T_A + T_B \circ \theta_{T_A}, \notag\\
T^*_A:= & \inf\{t > T_B: X_t \in A\} = T_B + T_A \circ \theta_{T_B}, \notag
\end{align}
and defining
\begin{align}
\psi_a := & \sum_{b' \in B} H_{B \cup C}(a,b') = P^a(T_B < T_C) \label{eq:psix}\\
\sigma_b := & \sum_{a' \in A} H_{A \cup C}(b,a') = P^b(T_A < T_C) \label{eq:sigmax}\\
\rho_a := & \sum_{b' \in B} H_{B \cup C}(a,b') \sigma_{b'} = P^a(T_B, T^*_A < T_C) \label{eq:rhox}\\
\phi_b := & \sum_{a' \in A} H_{A \cup C}(b,a') \psi_{a'} = P^b(T_A, T^*_B < T_C), \label{eq:phix}
\end{align}
we have the Green's function bounds 
\begin{align} 
G_{A}(a,a') \leq G_{A \cup B}(a,a') & \leq G_{A}(a,a') + \frac{\rho_a}{1-\rho_{a'}} G_{A}(a',a') \label{eq:GreenABUBa} \\
G_{B}(b,b') \leq G_{A \cup B}(b,b') & \leq G_{B}(b,b') + \frac{\phi_b}{1-\phi_{b'}} G_{B}(b',b') \label{eq:GreenABUBb} \\
0 \leq G_{A \cup B}(a,b) & \leq \min \left\{\frac{\sigma_b}{1 - \rho_a} G_{A}(a,a), \frac{\psi_a}{1-\phi_b} G_{B}(b,b) \right\}. \label{eq:GreenABUBab} 
\end{align}
\end{lem}
Recall that $G$ is symmetric, so the inputs can be swapped in any of these bounds. Also, by their definitions, $\psi_a \geq \rho_a$ for every $a \in A$ and $\sigma_b \geq \phi_b$ for every $b \in B$.

\pf See \cite[Proposition 1]{CarThree}. \qed

We now find the expected time of hitting the set $C$, starting from $A$, in terms of hitting $B \cup C$. Lower bounds are simple: just tack the other set on for a quicker hitting time. The upper bounds will require a recursive excursion treatment similar to the proof of Lemma \ref{lem:GreenABUB}.
\begin{lem} \label{lem:ExTCUB}
For $a \in A$ and $b \in B$, defining via \eqref{eq:psix} and \eqref{eq:sigmax}, 
\begin{equation} \label{eq:fAfB}
f_A := \sup_{a \in A} E^a(T_{B \cup C}), \,\,\, f_B := \sup_{b \in B} E^b(T_{A \cup C}), \,\,\, \psi := \sup_{a \in A} \psi_a, \,\,\, \sigma := \sup_{b \in B} \sigma_b,
\end{equation}
we have the expected hitting time bounds 
\begin{align}
E^a(T_{B \cup C}) \leq E^a(T_{C}) \leq E^a(T_{B \cup C}) + \psi_a\left[\frac{f_B  + \sigma f_A}{1 - \psi \sigma}\right] \label{eq:ExTCUBa} \\
E^b(T_{A \cup C}) \leq E^b(T_{C}) \leq E^b(T_{A \cup C}) + \sigma_b\left[\frac{f_A + \psi f_B}{1 - \psi \sigma}\right] \label{eq:ExTCUBb}
\end{align}
\end{lem}

\pf See \cite[Proposition 2]{CarThree}. \qed

\subsubsection{Application: Internal-External-Annulus Probabilities} \label{sec:4.1}

Let the following sets partition $\Z^2_K$, with $s \leq n < K \in \mathbb{N}$:
\begin{align*}
A = \hDn, \,\,\, B = \hDnsc, \,\,\, C = \hdDns.
\end{align*}

Starting from deep inside a disc, we first prove a bound on the probability of escaping the disc beyond an annulus outside it.
\begin{lem} \label{lem:BdEscapeEstimateToral}
\begin{align} 
\sup_{x \in D(0,n/2)} & P^x(T_{\bd D(0,n)_s} > T_{D(0,n+s)^c}) & \leq c(s^{-M+2} \lor n^{-M+2}). \label{eq:BdEscapeEstimatePlanar} \\
\psi = \sup_{\hat x \in \hp(D(0,n/2))} & P^{\hat x}(\hTdDns > T_{\hDnsc}) & \leq c (s^{-M+2} \lor n^{-M+2}). \label{eq:BdEscapeEstimateToral}
\end{align}
\end{lem}

\pf See \cite[Lemma 4.1]{CarVisits}. \qed

Note that for $\hat x \in \hDn$, by \eqref{eq:HittingTimeComp1},
\[ \{\hTdDns > \hTDnc\}^c = \{\hTdDns = \hTDnc\}.\]
Hence, provided $\hat x \in \hp(D(0,n/2))$, and $s \leq n$,  \eqref{eq:BdEscapeEstimateToral} is a bound for $\psi_{\hat x}$ from \eqref{eq:psix}.
Also, \eqref{eq:ProbZeroBeforeDisc} and \eqref{eq:BdEscapeEstimateToral} gives us the chance of escaping a disc, into its $s$-annulus, before visiting its center:
\begin{align} \label{eq:FP(2.48)}
& P^{\hat x}(T_{\hat{0}} > \hTDnc; \, \hTDnc = \hTdDns) \notag\\
 & = 1 - \frac{ \log(n/|\hat x|) + O(|\hat x|^{-1/4}) }{ \log n }(1 + O((\log n)^{-1})  + O(s^{-M+2}).
\end{align}

By \eqref{eq:BdEscapeEstimateToral}, \eqref{eq:rhox}, and \eqref{eq:phix}, for ${\hat x} \in \hp(D(0,n/2))$ and ${\hat y} \in \hDnsc$, 
\begin{align}
\rho_{\hat x} & = P^{\hat x}(T_{\hDnsc}, T_{\hDn}^* < T_{\hdDns}) \notag\\
& \leq c(s^{-M+2} \lor n^{-M+2}); \label{eq:rhoxBound} \\
\phi_{\hat y} & = P^{\hat y}(T_{\hDn}, T_{\hDnsc}^* < T_{\hdDns}) \notag\\
& \leq c(s^{-M+2} \lor n^{-M+2}). \label{eq:phixBound} 
\end{align}

Next, we find a 
bound for $\sigma_{\hat x}$ from \eqref{eq:sigmax}. 
\begin{lem} \label{lem:FPLemma2.6}
For $n$ sufficiently large, 
\begin{align}
\sigma & = \sup_{\hat x \in \hDnsc} P^{\hat x}(T_{\hDn} < T_{\hdDns}) \notag \\
 & \leq cn^2 \log(n)^2 (s^{-M} + n^{-M}). \label{eq:FP(2.71)}
\end{align}
\end{lem}

\pf See \cite[Lemma 4.2]{CarVisits}. \qed

In particular, if $s=O(n)$, since $M=4+2\beta$, \eqref{eq:FP(2.71)} is bounded above by $cn^{-2}$, and if $s=O(\sqrt{n})$, \eqref{eq:FP(2.71)} is bounded above by $cn^{-\beta}$.

Combining \eqref{eq:ToralGamblersSuccess}  and \eqref{eq:BdEscapeEstimateToral}, we find the probability that, starting far from a small disc $\hDr$, the walk escapes a larger disc $\hDR$ before entering $\hDr$. If $r < R$ and $\hat x \in \hp(D(0,R/2))$,
we have 
\begin{align} \label{eq:FP2.49analog}
& P^{\hat x}(\hTDRc < \hTDr; \hTDRc = T_{\hp(\bd D(0,R)_s)})\notag\\
& = \frac{\log(|\hat x|/r) + O(r^{-1/4})}{\log(R/r)} + O(s^{-M+2}).
\end{align}

To enter a disc, we first quote the planar result \cite[Lemma 2.4]{BRFreq}: if $s < r < R$ sufficiently large with $R \leq r^2$ we can find $c < \infty$ and $\delta > 0$ such that for any $ r < |x| < R$, 
\begin{align} \label{eq:FP(2.50)}
P^x(T_{D(0,r)} < T_{D(0,R)^c}; T_{D(0,r)} = T_{D(0,r-s)}) \leq cr^{-\delta} + cs^{-M+2}.
\end{align}
We see the same result on $\Z^2_K$, with an extra toral term (which is absorbed). 
\begin{lem} \label{lem:FPLemma2.4}
For the conditions listed above, 
\begin{align} 
P^{\hat x}(\hTDr & < \hTDRc; \hTDr = T_{\hp(D(0,r-s))}) \notag \\
 & \leq cr^{-\delta} + cs^{-M+2}. \label{eq:FP2.50analog}
\end{align}
\end{lem}

\pf See \cite[Lemma 4.3]{CarVisits}. \qed

We use \eqref{eq:FP2.50analog} along with \eqref{eq:ToralGamblersRuin} to get the toral gambler's ruin-via-annulus estimate:  
\begin{align} \label{eq:FP2.66analog}
& P^{\hat x}(\hTDr < \hTDRc; \hTDr = T_{\hp(\bd D(0,r-s)_s)}) \notag\\
& = \frac{\log(R/|\hat x|) + O(r^{-\delta})}{\log(R/r)} + O(s^{-M+2}).  
\end{align}

We now give results on these probabilities for a finely-tuned set of radii and annuli which will appear in later sections. For $n$ large and $c > 0$ and set the following:\footnote{The use of different thicknesses of $s_{n-1}$ depending on direction is due to the entry probability from level $n$ in the lower bound argument of Section \ref{ch:Late}; see Section \ref{sec:xband} and \eqref{eq:bn0} for details.}
 \index{rnk@$r_{n,k}$} 
\begin{align}
r_{n,k} = e^n n^{3k}, s_k = n^4, & & r_{n,k}' = r_{n,k}+s_{k},  & & k=0,1, \ldots, n; \notag\\
s_{n-1}^{n \downarrow} = \sqrt{r_{n,n-1}}. \label{eq:rsKdefn}
\end{align}
For large enough $n$, $n^4 < r_{n,l}^{\delta}$ for any $1/2 \leq \delta < 1$, so for any $\hat x \in \hp(\bd D(0,r_{n,l})_{s_l})$ and $1 \leq l \leq n-1$, 
\begin{align}
 n^3 = \frac{r_{n,l}}{r_{n,l-1}} & < \frac{|\hat x|}{r_{n,l-1}} < \frac{r_{n,l}+r_{n,l}^{\delta}}{r_{n,l-1}} < n^3 + e^{-n(1-\delta)} n^{-3l(1-\delta)+3} < n^3 + n^{-1} \notag\\
\implies & \log\left(\frac{|\hat x|}{r_{n,l-1}}\right) = 3 \log n + O(n^{-4}), \label{eq:xbandlogbound}
\end{align}
so by \eqref{eq:FP2.49analog} and \eqref{eq:xbandlogbound} we have 
\begin{align} 
a_{l+1} := & \, P^{\hat x}\left( T_{\hp(D(0,r_{n,l+1})^c)} < T_{\hp(D(0,r_{n,l-1}'))}; T_{\hp(D(0,r_{n,l+1})^c)} = T_{\hp(\bd D(0,r_{n,l+1})_{s_{l+1}})} \right)\notag\\
 = & \, \frac{3 \log n + O(n^{-4}) + O(r_{n,l-1}^{-1/4})}{\log(r_{n,l+1}/r_{n,l-1})} + O(s_l^{-M+2}) \label{eq:FP6.7analog} \\ 
= & \, \frac{3 \log n + O(n^{-4})}{6 \log n} + O(s_l^{-M+2}) = \frac{1}{2} + o(n^{-4}), \notag
\end{align}
Likewise, using \eqref{eq:xbandlogbound}, 
\begin{align}
 n^{-3} = \frac{r_{n,l}}{r_{n,l+1}} & < \frac{|\hat x|}{r_{n,l+1}} < \frac{r_{n,l}+r_{n,l}^{\delta}}{r_{n,l+1}} < n^{-6}(n^3 + e^{-n(1-\delta)} n^{-3l(1-\delta)+3}) < n^{-3} + n^{-7} \notag\\
\implies & n^3 - \frac{n^3}{n^4+1} = n^3 - O(n^{-1}) < \frac{r_{n,l+1}}{|\hat x|} < n^3 \label{eq:xbandlogbound2}\\
\implies & \log\left(\frac{r_{n,l+1}}{|\hat x|}\right) = 3 \log n + O(n^{-4}), \notag
\end{align}
so by \eqref{eq:FP2.66analog} and \eqref{eq:xbandlogbound2} we have  
\begin{align}
b_{l} := & \, P^{\hat x}\left(T_{\hp(D(0,r_{n,l-1}'))} < T_{\hp(D(0,r_{n,l+1})^c)}; T_{\hp(D(0,r_{n,l-1}'))} = T_{\hp(\bd D(0,r_{n,l-1})_{s_{l-1}})} \right) \notag\\
 = & \, \frac{1}{2} + o(n^{-4}). \label{eq:FP6.8analog}
\end{align}

\subsubsection{Application: Green's Functions, Hitting Times}

We start calculating bounds for the external Green's function with $\hat{x} \in \hp(D(0,n/2))$, $\hat{y} \in \hDn$: by \eqref{eq:GreenABUBa} with $A=\hDn$, \eqref{eqn:GreenXZBounds}, and \eqref{eq:rhoxBound}, 
\begin{align}
\hat{G}_{\hdDnsc}(\hat x,\hat y) & \leq \hat{G}_{\hDn}(\hat x, \hat y) + \frac{\rho_{\hat x}}{1 - \rho_{\hat y}} \hat{G}_{\hDn}(\hat y, \hat y). \label{eq:GhDnBound}  
\end{align}
In particular, if $\hat{y} = \hat{0}$ and $s = O(n)$, then $\rho_{\hat x} \leq c n^{-2}$ and by \eqref{eq:FP(2.18)}, 
\begin{align}
\hat{G}_{\hdDnsc}(\hat x,\hat 0) & \leq \hat{G}_{\hDn}(\hat x, \hat 0) + \frac{\rho_{\hat x}}{1 - \rho_{\hat 0}} \hat{G}_{\hDn}(\hat 0, \hat 0) \notag\\
\implies \hat{G}_{\hdDnsc}(\hat x,\hat 0) & = \frac{2}{\pi_{\Gamma}}\log\left(\frac{n}{|\hat x|}\right)  + C(\hat p_1) + O(|\hat x|^{-1/4}). \label{eq:GhDnBoundZero} 
\end{align}
By \eqref{eq:GreenABUBb}, \eqref{eqn:ExternalGreenToralComplete}, and \eqref{eq:phixBound}, for $\hat{x}, \hat{y} \in \hDnsc$, 
\begin{align}
\hat{G}_{\hdDnsc}(\hat x,\hat y) & \leq \hat{G}_{\hDnsc}(\hat x, \hat y)  + \frac{\phi_{\hat x}}{1 - \phi_{\hat y}} \hat{G}_{\hDnsc}(\hat y, \hat y)  \notag\\
 & \leq c (\log(|\hat{x}| \land |\hat{y}|))^2. \label{eq:GhDnscBound} 
\end{align}
Finally, for $\hat{x} \in \hp(D(0,n/2))$ and $\hat{y} \in \hDnsc$, by \eqref{eq:GreenABUBab}, \eqref{eq:BdEscapeEstimateToral}, \eqref{eq:FP(2.71)}, and the above, 
\begin{align}
\hat{G}_{\hdDnsc}&(\hat x,\hat y) \label{eq:GreenbdDUB}\\
\leq & \min\left\{ \frac{\sigma_x}{1 - \rho_y} \hat{G}_{\hDn}(\hat x, \hat x), 
 \frac{\psi_x}{1-\phi_y} \hat{G}_{\hDnsc}(\hat y, \hat y) \right\} \notag\\
 \leq & \, c \, \min\left\{ n^2 (\log n)^3 (s^{-M} + n^{-M}), (\log(|\hat{y}|))^2  (s^{-M+2} \vee n^{-M+2}) \right\}. \notag
\end{align}
In particular, if $s=O(n)$, then in this case $\hat{G}_{\hdDnsc}(\hat x,\hat y) \leq c n^{-2}$, and if $s=O(\sqrt{n})$, the bound is $cn^{-\beta}$. 

By \eqref{eq:HittingTimeComp2} and \eqref{eqn:EnterDiscExp}, for $y \in D(0,n+s)^c \subset \Z^2$, the external planar annulus hitting time $E^y(T_{\bd D(0,n)_s}) = \infty$. Since, starting from inside the disc $x \in D(0,n)$, there is positive probability of hopping over an $s$-width annulus, then by the strong Markov property on $T_{D(0,n+s)^c}$, the internal planar annulus hitting time $E^x(T_{\bd D(0,n)_s}) = \infty$ as well. This is not the case for the toral analogues of these times. 

Torally, our walk can make small or targeted jumps before the disc escape time. To bound the annulus hitting times, we  employ \eqref{eqn:EscapeDiscExpTorus}, \eqref{eq:ToralDiscEntryExpectedTime}, and \eqref{eq:fAfB}. These yield, for some $c,c' < \infty$, 
\begin{align}
f_{\hDn}     = & \sup_{\hat x \in \hDn} E^{\hat x}(T_{\hDnc}) \leq cn^2, \label{eq:fD0n}\\
f_{\hDnsc} = & \sup_{\hat y \in \hDnsc} E^{\hat y}(T_{\hDns}) \leq c' (K \log K)^2. \label{eq:fD0nsc}
\end{align}
By \eqref{eq:ExTCUBa}, \eqref{eq:ExTCUBb}, \eqref{eq:fD0n}, \eqref{eq:fD0nsc}, \eqref{eq:BdEscapeEstimateToral}, and \eqref{eq:FP(2.71)}, the expected annulus hitting time is bounded above: if $\hat{x} \in \hp(D(0,n/2))$ and $\hat{y} \in \hDnsc$, 
\begin{align} 
\EV^{\hat x}(T_{\hdDns}) & \leq \EV^{\hat x}(T_{\hDnc}) + \psi_{\hat x} \left[ \frac{f_{\hDnsc} + \sigma f_{\hDn}}{1 - \psi\sigma} \right] \notag\\
 & \leq \EV^{\hat x}(T_{\hDnc}) + c\left(s^{-M+2} \vee n^{-M+2}\right)(K \log K)^2; \label{eq:InTbdDnsBounds}\\
\EV^{\hat y}(T_{\hdDns}) & \leq \EV^{\hat y}(T_{\hDns}) + \sigma_{\hat y} \left[ \frac{f_{\hDn} + \psi f_{\hDnsc}}{1 - \psi\sigma} \right] \notag\\
 & \leq c(K \log K)^2. \label{eq:ExTbdDnsBounds}
\end{align}
In particular, if $s, n=O(K)$, then for $K$ sufficiently large, note that by \eqref{eqn:EscapeDiscExpTorus}, 
\begin{align*}
\EV^{\hat x}(T_{\hDnc}) = \frac{K^2 - |\hat x|^2}{\gamma^2} + O(K), 
\end{align*}
which, with $M = 4 + 2\beta$, reduces \eqref{eq:InTbdDnsBounds} to 
\begin{align} 
\EV^{\hat x}(T_{\hdDns}) & = \left(1 + O(K^{-2-\beta})\right) \EV^{\hat x}(T_{\hDnc}). \label{eq:InTbdDnsBounds2}
\end{align}

\section{Harnack Inequalities} \label{ch:Harnack}

Here we will quote and apply Harnack inequality results from \cite{CarHarnack} for use in our excursion treatments.

\subsection{Interior Harnack inequalities}

Our first interior Harnack inequality gives estimates on the probability, when escaping a large disc from deep inside it, of landing in an annulus close to the disc's boundary. 

\begin{prop} \label{lem:InteriorHarnack2}
Uniformly for $1 \leq m \ll r$, with $s \ll \frac{r}{4m}$, $x, x' \in D(0,2r)$, $R=4mr$, and $y \in D(0,R)^c$,   
\begin{align} 
H_{D(0,R)^c}(x,y) & = (1 + O(m^{-1})) H_{D(0,R)^c}(x',y) + O(R^{-M} \log R), \label{eq:IntHarn2}
\end{align}
where the error term is completely absorbed, \ie 
\begin{align} 
H_{D(0,R)^c}(x,y) & = (1 + O(m^{-1})) H_{D(0,R)^c}(x',y), \label{eq:IntHarn2NoError}
\end{align}
if $s \leq (\log R)^4$ and $y \in \bd D(0,R)_{s}$. 

Furthermore, if $x \in \bd D(0,r)_{r}$ and $y \in D(0,R)^c$, 
\begin{align} 
 & P^x \left(S_{T_{D(0,R)^c}}=y, \, T_{D(0,R)^c} < T_{D\left(0,\frac{r}{4m}+s\right)}\right) \label{eq:IntHarn2OutFirst} \\
  & \,\,\, = (1 + O(m^{-1})) P^x \left(T_{D(0,R)^c} < T_{D\left(0,\frac{r}{4m}+s\right)}\right) H_{D(0,R)^c}(x,y) + O(R^{-M} \log R), \notag 
\end{align}  
with a similar loss of the error term if $y \in \bd D(0,R)_{s}$. 
\end{prop}

\pf See \cite[Prop. 3.1]{CarHarnack}. \qed

Here is a focused result for our applications which follows directly.

\begin{cor} \label{lem:InteriorHarnack}
Let $e^n \leq r$, $R = n^3 r$ (\ie $m=\frac{n^3}{4}$ for $R = 4mr$). Uniformly for $x, x' \in D(0,r+\sqrt{r})$ and $y \in \bd D(0,R)_{n^4}$,  
\begin{equation} \label{eq:InteriorHarnackPlanar}
H_{D(0,R)^c}(x,y) = \left(1 + O\left(n^{-3}\right)\right) H_{D(0,R)^c}(x',y).
\end{equation}
Furthermore, uniformly in $x \in \bd D(0,r)_{\sqrt{r}}$ and $y \in \bd D(0,R)_{n^4}$, 
\begin{align}
& P^x(S_{T_{D(0,R)^c}} = y, T_{D(0,R)^c} < T_{D(0,\frac{r}{n^3}+n^4)}) \label{eq:InteriorHarnackPlanarZoomIn}\\
& \, = \left(1 + O\left(n^{-3}\right)\right) P^x(T_{D(0,R)^c} < T_{D(0,\frac{r}{n^3}+n^4)}) H_{D(0,R)^c}(x,y). \notag
\end{align}
\end{cor}

We now move these results to the torus.

\begin{prop} \label{lem:InteriorHarnack2Toral}
For large $r$ and $1 \leq m \ll r$ such that $R = 4mr < K/6$ and $s \leq (\log R)^4$, uniformly for $\hat x, \hat x' \in \hp(D(0,2r))$ and $\hat y \in \hDRc$, 
\begin{align} 
\hat{H}_{\hDRc}(\hat x,\hat y)
 = & \left(1 + O\left(m^{-1}\right)\right) \hat{H}_{\hDRc}(\hat x',\hat y) \notag\\
 & + O(R^{-M} \log R \lor K^{-M} R^2). \label{eq:IntHarn2Toral}
\end{align}
Furthermore, uniformly in $\hat x \in \hp(\bd D(0,r)_{r})$ and $\hat y \in \hDRc$,
\begin{align}
 P^{\hat x}&(\hat{S}_{\hTDRc} = \hat y, \hTDRc < T_{\hp(D(0,\frac{r}{4m}+s))}) \notag\\
 & = \left(1 + O\left(m^{-1}\right)\right) P^{\hat x}(\hTDRc < T_{\hp(D(0,\frac{r}{4m}+s))}) \hat{H}_{\hDRc}(\hat x,\hat y) \notag\\
 & \,\, + O(R^{-M} \log R \lor K^{-M} R^2). \label{eq:IntHarn2ToralZoomIn}
\end{align}

If $\hat y \in \hp(\bd D(0,R)_{s})$, 
the error term is absorbed in both of these statements. 
\end{prop}

\pf See \cite[Prop. 3.2]{CarHarnack}. \qed

\begin{cor} \label{lem:InteriorHarnackToral}
Let $n > 13$, $e^n \leq r$, $R = n^3 r$ (\ie $m=\frac{n^3}{4}$ for $R = 4mr$). Uniformly for $\hat x, \hat x' \in \hp(D(0,2r))$, $K > 4(R + n^4)$, and $\hat y \in \hp(\bd D(0,R)_{n^4})$, 
\begin{equation} \label{eq:InteriorHarnack}
\hat{H}_{\hDRc}(\hat x,\hat y) = \left(1 + O\left(n^{-3}\right)\right) \hat{H}_{\hDRc}(\hat x',\hat y).
\end{equation}
Furthermore, uniformly in $\hat x \in \hp(\bd D(0,r)_{\sqrt{r}})$ and $\hat y \in \hp(\bd D(0,R)_{n^4})$, 
\begin{align} 
P^{\hat x}&(\hat{S}_{\hTDRc} = \hat y, \hTDRc < T_{\hp(D(0,\frac{r}{n^3}+n^4))}) \label{eq:InteriorHarnackZoomIn}\\
= & \left(1 + O\left(n^{-3}\right)\right) P^{\hat x}(\hTDRc < T_{\hp(D(0,\frac{r}{n^3}+n^4))} ) \hat{H}_{\hDRc}(\hat x,\hat y). \notag
\end{align}
\end{cor}

\subsection{Exterior Harnack inequality}

 We now give general and applied Harnack inequalities for the plane and torus dealing with entering a small disc from far outside. 

\begin{prop} \label{lem:ExteriorHarnackPlanarGeneral}
Let $R = 4mr$ with $1 \leq m \ll r$ ($m = o(r^{1/4})$) and large enough $r$, and $s \leq (\log R)^4$. Then, uniformly for $x, x' \in D(0,R)^c$ and $y \in \bd D(0,r)_{s}$, 
\begin{align} 
H_{D(0,r+s)}(x,y) & = \left(1 + O\left(m^{-1} \log m\right)\right) H_{D(0,r+s)}(x',y). \label{eq:ExteriorHarnackPlanarGeneral}
\end{align}
Furthermore, for $x, x' \in \bd D(0,R)_{\sqrt{R}}$, 
\begin{align}
P^x& (S_{T_{D(0,r+s)}} = y; \, T_{D(0,r+s)} < T_{{D(0,4mR)}^c}) \label{eq:ExteriorHarnackPlanarZoomOutGeneral} \\
&  = \left(1 + O\left(m^{-1} \log m\right)\right) H_{D(0,r+s)}(x,y) P^x(T_{D(0,r+s)} < T_{{D(0,4mR)}^c})  \notag\\
&  = \left(1 + O\left(m^{-1} \log m\right)\right) P^{x'}(S_{T_{D(0,r+s)}} = y; T_{D(0,r+s)} < T_{{D(0,4mR)}^c}). \notag
\end{align}
\end{prop}

\pf See \cite[Prop. 4.1]{CarHarnack}. \qed

We now fine-tune this result for our applications

\begin{cor} \label{lem:ExteriorHarnack}
As in Lemma \ref{lem:InteriorHarnack}, let $e^n \leq r$, $R = 4mr = n^3 r$. Then, uniformly for $x, x' \in D(0,R)^c$ and $y \in \bd D(0,r)_{n^4}$, 
\begin{align} 
H_{D(0,r+n^4)}(x,y) & = \left(1 + O\left(n^{-3} \log n\right)\right) H_{D(0,r+n^4)}(x',y). \label{eq:ExteriorHarnackPlanar}
\end{align}
Furthermore, for $x, x' \in \bd D(0,R)_{\sqrt{R}}$, 
\begin{align}
P^x(S_{T_{D(0,r+n^4)}} = y; \, & T_{D(0,r+n^4)} < T_{{D(0,n^3 R)}^c}) \label{eq:ExteriorHarnackPlanarZoomOut} \\
&  = \left(1 + O\left(n^{-3} \log n\right)\right) H_{D(0,r+n^4)}(x,y) P^x(T_{D(0,r+n^4)} < T_{{D(0,n^3 R)}^c}) \notag\\
&  = \left(1 + O\left(n^{-3} \log n\right)\right) P^{x'}(S_{T_{D(0,r+n^4)}} = y; T_{D(0,r+n^4)} < T_{{D(0,n^3 R)}^c}). \notag
\end{align}
\end{cor}

\newcommand{\Drs}{D(0,r+s)}
\newcommand{\hDrs}{\hp(D(0,r+s))}
\newcommand{\TDrs}{T_{D(0,r+s)}}
\newcommand{\hTDrs}{T_{\hp(D(0,r+s))}}

When attempting to move the planar exterior Harnack inequality to the torus, we run into difficulties in dealing with walks that wander and enter far-off copies of $D(0,r+s)$ instead of the primary copy. We modify the exterior Harnack inequality for the toral case to fit our requirements; \eqref{eq:ExteriorHarnackToralZoomOut} is a direct application of \eqref{eq:ExteriorHarnackToralZoomOutGeneral}. 

\begin{prop} \label{lem:ExteriorHarnackToral}
Let $R = 4mr$ with $1 \leq m = o(r^{1/4})$ and large enough $r$, $4mR < K/4$, and $s \leq (\log R)^4$. Then, uniformly for $\hat x, \hat x' \in \hp(\bd D(0,R)_{\sqrt{R}})$ and $\hat y \in \hp(\bd D(0,r)_{s})$, 
\begin{align}
P^{\hat x}& (\hat{S}_{\hTDrs} = \hat y; \, \hTDrs < T_{{\hp(D(0,4mR)}^c_K)}) \label{eq:ExteriorHarnackToralZoomOutGeneral} \\
&  = \left(1 + O\left(m^{-1} \log m\right)\right) P^{\hat x'}(\hat{S}_{T_{\hp(D(0,r+s))}} = \hat y; T_{\hp(D(0,r+s))} < T_{{\hp(D(0,4mR)}^c_K)}). \notag
\end{align}
\end{prop}

\pf See \cite[Prop. 4.2]{CarHarnack}. \qed

\begin{cor} \label{lem:ExteriorHarnackToralZoomOut}
Let $e^n \leq r$, $R = 4mr = n^3 r$. Then, uniformly for $\hat x, \hat x' \in \hp(\bd D(0,R)_{\sqrt{R}})$ and $\hat y \in \hp(\bd D(0,r)_{n^4})$, 
\begin{align}
P^{\hat x}&(\hat{S}_{T_{\hp(D(0,r+n^4))}} = \hat y; \, T_{\hp(D(0,r+n^4))} < T_{{\hp(D(0,n^3 R)}^c_K)}) \label{eq:ExteriorHarnackToralZoomOut} \\
&  = \left(1 + O\left(n^{-3} \log n\right)\right) P^{\hat x'}(\hat{S}_{T_{\hp(D(0,r+n^4))}} = \hat y; T_{\hp(D(0,r+n^4))} < T_{{\hp(D(0,n^3 R)}^c_K)}). \notag
\end{align}
\end{cor}


\section{Excursions} \label{ch:Excursions}

In this section we find bounds on times of excursions between concentric annuli. As in \cite{DPRZ2006}, for any hitting time $\hat T$ on the torus $\Z^2_K$, we set
\[||\hat T|| := \sup_{\hat y \in \Z^2_K} \EV^{\hat y}(\hat T).\]
By Kac's moment formula for the strong Markov process $\hat{S}_t$ (see \cite[(6)]{FitzPit}), we have for any $t$ and $\hat y$,
\begin{equation} \label{eq:KacMomentStrongMarkov}
\EV^{\hat y}(\hat T^k) \leq k! \EV^{\hat y}(\hat T) ||\hat T||^{k-1}.
\end{equation}

\subsection{Between a small annulus and far out}

Let $R = 4mr$. In this section, when considering visits to $\hat x \in \Z^2_K$, we will consider excursions between a small annulus and the complement of a large disc, both centered at $\hat x$. Define the times 
\begin{eqnarray}
\tau^{(0)} & = & \inf \{t \geq 0: \hat{S}_t \in \Dr\},\\
\sigma^{(1)} & = & \inf \{t \geq \tau^{(0)}: \hat{S}_t \in \DRc\},
\end{eqnarray}
and inductively for $j=1,2,\ldots$, let 
\begin{eqnarray}
\tau^{(j)} & = & \inf \{t \geq \sigma^{(j)}: \hat{S}_{t+\mathfrak{T}_{j-1}} \in \Dr\},\\
\sigma^{(j+1)} & = & \inf \{t \geq 0: \hat{S}_{t+\mathfrak{T}_{j}} \in \DRc\},
\end{eqnarray}
where $\mathfrak{T}_j = \sum_{i=0}^j \tau^{(i)}$ for $j=0,1,2,\ldots$. Thus $\tau^{(j)}$ is the length of time of the $j$th excursion $\E_j$ from $\Dr \rightarrow \DRc \rightarrow \Dr$, and $\sigma^{(j)}$ is the amount of time it takes for the first leg of $\E_j$. From here on, set $\tau = \tau^{(1)}$.

\begin{figure}[!ht]
  \centering
    \includegraphics[width=3in]{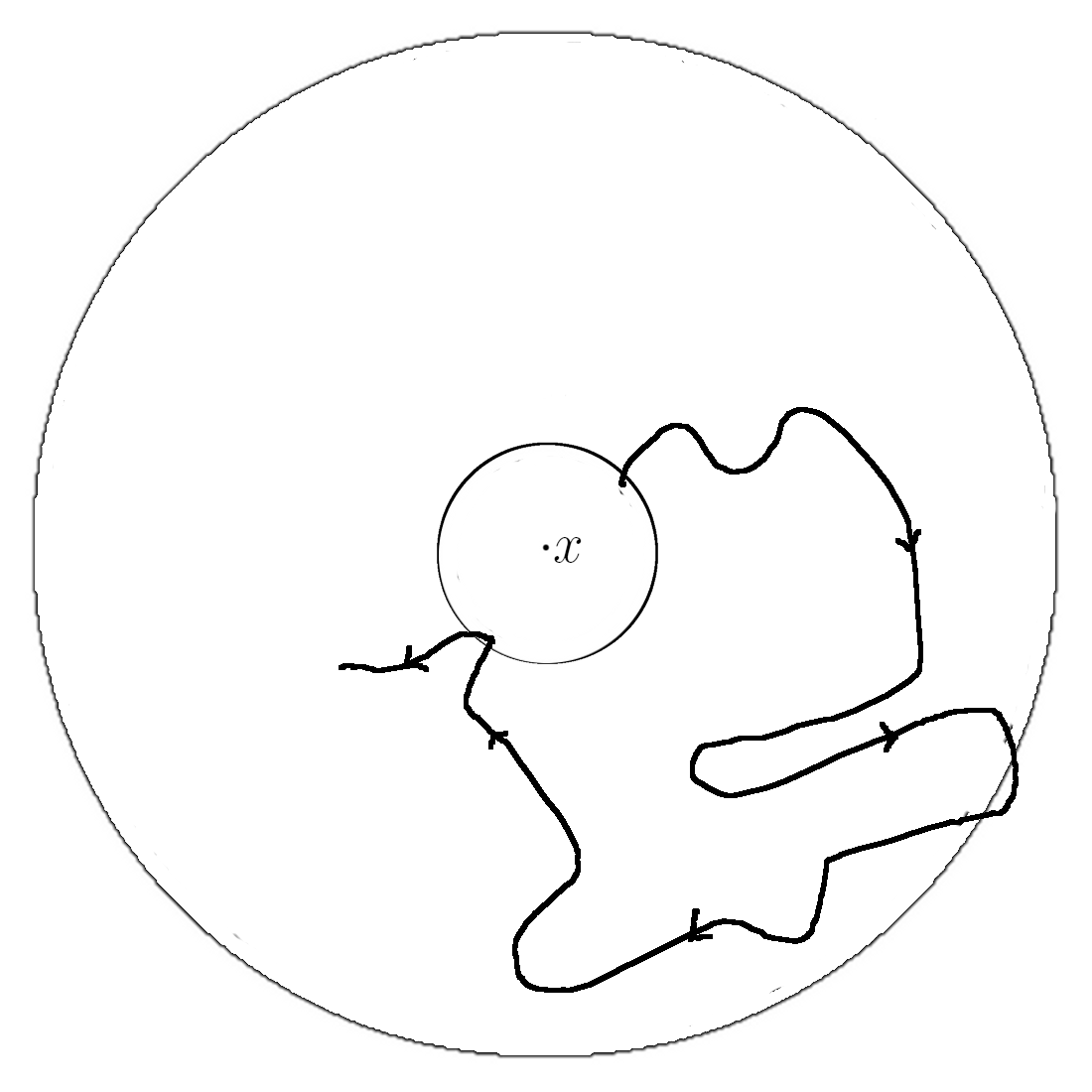}
    \parbox{4in}{
    \caption{A sample excursion $\E_j$.}
    \label{fig:rstar_Rstar_excursion}}
\end{figure}

Our first lemma gives bounds on these excursion times, and shows their concentration near the asymptotic limit. 

\begin{lem} \label{lem:LPLemma3.1}
Uniformly for $1 \leq m < r$, $R = 4mr$, 
$cK^{1-\epsilon} = R \leq \frac{K}{24}$ for some small $0 \leq \epsilon \ll \min\{\beta, \frac{1}{2}\}$, 
and $(\log K)^2 < s < (\log R)^4$, $\exists \, c_1 < \infty$ \st $\forall \eta$: $1 \geq \eta \geq c_1 \left(\left(\frac{r}{R}\right) + s^{-1} + K^{-2\beta-2\epsilon}(\log K)^2\right)$, 
\begin{align} \label{eq:ExcursionLemma1}
(1 - \eta) \frac{2}{\pi_{\Gamma}} K^2 \log\left(\frac{R}{r}\right) & \leq \min_{\hat x,\hat y \in \Z^2_K} \EV^{\hat y}(\tau) \\
& \leq \max_{\hat x,\hat y \in \Z^2_K} \EV^{\hat y}(\tau) \leq (1 + \eta) \frac{2}{\pi_{\Gamma}} K^2 \log\left(\frac{R}{r}\right). \notag
\end{align}
\end{lem}

\pf Note that $\hat x$ is the center of the discs we will analyze. 
Let $\hat{S}_0$ be distributed uniformly on $\ZK$. Then $\{\hat{S}_t\}$ is a stationary and ergodic stochastic process. By Birkhoff's ergodic theorem we then have that 
\[\lim_{t \to \infty} {1 \over t} \sum_{i=0}^t 1_{\{\hat x\}}(\hat{S}_i) = {1 \over K^2}
\,\,\,\,\,\, \mbox{ a.s.}\]
Thus, with $\mathfrak{T}_{-1} = 0$,
\begin{equation} \label{eq:ExcursionLemma1Limit1}
\lim_{t \to \infty} \frac{{1 \over t} \sum_{j=0}^t \sum_{i=0}^{\tau^{(j)}}
1_{\{\hat x\}}(\hat{S}_{i+\mathfrak{T}_{j-1}})}{{1 \over t} \sum_{j=0}^t \tau^{(j)}} =
\frac{1}{K^2}  \,\,\,\,\,\, \mbox{ a.s.}
\end{equation}
Let $\rho$ be uniform measure on $\ZK$, and for $j \geq 1$, let
\[Z_j := \tau^{(j)} - \EV^{\rho}(\tau^{(j)}|\F_{\mathfrak{T}_{j-1}}) = \tau^{(j)} - \EV^{\hat{S}_{\mathfrak{T}_{j-1}}}(\tau).\]
By the strong Markov property, $\{Z_j\}$ is an orthogonal sequence. Since any irreducible, aperiodic Markov chain with finite state space is positive recurrent, we have that $||T_{\Dr}||$, $||T_{\DRc}|| < \infty$, and using \eqref{eq:KacMomentStrongMarkov} we see that the sequence $\{\tau^{(j)}\}$ and hence $\{Z_j\}$ has uniformly bounded second moments. It follows from Rajchman's strong law of large numbers that
\begin{equation} \label{eq:ExcursionLemma1Sublimittau}
\lim_{t \to \infty} {1 \over t} \sum_{j=1}^t [\tau^{(j)} - \EV^{\hat{S}_{\mathfrak{T}_{j-1}}}(\tau)]  = 0 \,\,\,\,\,\, \mbox{ a.s.}
\end{equation}
Similarly, set $\sigma^{(0)} = \tau^{(0)}$ and for $j \geq 0$ let $Y_j$ be the number of visits to $\hat x$ on the $j$th excursion $\Dr \to \DRc \to \Dr$: 
\begin{equation} \label{eq:ExcursionDecomp}
Y_j := \sum_{i=0}^{\tau^{(j)}} 1_{\{\hat x\}}(\hat{S}_{i + \mathfrak{T}_{j-1}})
  = \sum_{i=0}^{\sigma^{(j)}} 1_{\{\hat x\}}(\hat{S}_{i + \mathfrak{T}_{j-1}})
 + \sum_{i=\sigma^{(j)}+1}^{\tau^{(j)}} 1_{\{\hat x\}}(\hat{S}_{i + \mathfrak{T}_{j-1}}).
\end{equation}
Define
\[\tilde{Y}_j := Y_j - \EV^{\rho}(Y_j|\F_{\mathfrak{T}_{j-1}}) = Y_j - \EV^{\hat{S}_{\mathfrak{T}_{j-1}}}(Y_1).\]
By the strong Markov property, $\{\tilde{Y}_j\}$ is also an orthogonal sequence, and since $Y_j \leq \tau^{(j)}$, the sequence $\{\tilde{Y}_j\}$ also has uniformly bounded second moments. Thus, by Rajchman's strong law of large numbers, 
\begin{equation} \label{eq:ExcursionLemma1SublimitY}
\lim_{t \to \infty} {1 \over t} \sum_{j=1}^t [Y_j - \EV^{\hat{S}_{\mathfrak{T}_{j-1}}}(Y_1)]  = 0 \,\,\,\,\,\, \mbox{ a.s.}
\end{equation}
Let $\hat y \in \Dr$. To bound $\EV^{\hat y}(Y_1)$ we need to consider the two sums in \eqref{eq:ExcursionDecomp}. By \eqref{eqn:Green}, \eqref{eq:ExcursionDecomp}, and the strong Markov property at $\sigma^{(1)}$, we have 
\[\EV^{\hat y}(Y_1) = \hat{G}_{\DRcc}(\hat y,\hat x) + \EV^{\hat y}\left(\hat{G}_{\Drc}\left(\hat{S}_{T_{\DRc}},\hat x\right)\right).\]
By \eqref{eq:GreenXZeroTorusCalc}, for some constant $c^* = c^*(\hat{p}_1)$, and any $\hat y \in \hp(\bd D(x,r)_s)$, 
\[\hat{G}_{\DRcc}(\hat y,\hat x) = \frac{2}{\pi_{\Gamma}} \log \left( \frac{R}{r} \right) + c^* + O(r^{-1/4}).\] 
Also, 
$O(R) \leq |\hat{S}_{T_{\DRc}} - \hat x| \leq O(K)$, so by \eqref{eq:GreenbdDUB} and $(\log K)^2 < s$, 
\begin{align*}
\EV^{\hat y}\left(\hat{G}_{\Drc}\left(\hat{S}_{T_{\DRc}}, \hat x\right)\right) & \leq c (\log K)^2 s^{-M+2} \leq c s^{-M+3} = o(s^{-1}).
\end{align*}
Hence, for some finite universal constant $c_0>0$ and all allowable $s$,
\begin{align} 
\frac{2}{\pi_{\Gamma}} \log\left(\frac{R}{r}\right) + c^* - c_0 s^{-1} & \leq \min_{\hat x} \min_{\hat y \in \Dr} \EV^{\hat y}(Y_1) \label{eq:ExcursionBoundsMinMinMaxMax}\\
  & \leq \max_{\hat x} \max_{\hat y \in \Dr} \EV^{\hat y}(Y_1) \leq \frac{2}{\pi_{\Gamma}} \log\left(\frac{R}{r}\right) + c^* + c_0 s^{-1}.  \notag
\end{align}
With $\tau^{(0)}$ finite, we get by combining (\ref{eq:ExcursionLemma1Limit1}), (\ref{eq:ExcursionLemma1Sublimittau}), and (\ref{eq:ExcursionLemma1SublimitY}) that, a.s.,
\begin{equation} \label{eq:ExcursionTimeToVisitsRatio}
\lim_{t \to \infty} \frac{\frac{1}{t} \sum_{j=1}^t \EV^{\hat{S}_{\mathfrak{T}_{j-1}}}(\tau)}
 {\frac{1}{t} \sum_{j=1}^t \EV^{\hat{S}_{\mathfrak{T}_{j-1}}}(Y_1)} = K^2.
\end{equation}
Consequently, in view of \eqref{eq:ExcursionBoundsMinMinMaxMax}, for some universal constant $c_2$ and all 
$1 \geq \eta \geq c_2 \left(s^{-1} + \frac{r}{R}\right)$,  
\begin{align} \label{eq:OuterExcursionBounds}
\min_{\hat y \in \Dr} \EV^{\hat y}(\tau)  & \leq \frac{2}{\pi_{\Gamma}} K^2 \left(1 + \frac{\eta}{3}\right)\log\left(\frac{R}{r}\right) \notag\\
\max_{\hat y \in \Dr} \EV^{\hat y}(\tau) & \geq \frac{2}{\pi_{\Gamma}} K^2 \left(1 - \frac{\eta}{3}\right) \log\left(\frac{R}{r}\right)
\end{align}
For $\hat y \in \Dr$ we have $\tau^{(0)} = 0$ and by the strong Markov property at $\sigma^{(1)}$,
\begin{equation} \label{eq:ExpValTau}
\EV^{\hat y}(\tau) = \EV^{\hat y}(T_{\DRc}) + \sum_{\hat z \in \DRc} \hat{H}_{\DRc}(\hat y,\hat z) \, \EV^{\hat z}(T_{\Dr}).
\end{equation}
By \eqref{eqn:EscapeDiscExpTorus} and $R = cK^{1-\epsilon}$, 
\begin{align}
\EV^{\hat y}(T_{\DRc}) = cK^{2-2\epsilon} + O(K^{1-\epsilon}) \label{eq:DRcBounds}
\end{align}
for every $\hat y \in \Dr$. Hence, 
\begin{align}   
\max_{\hat y \in \Dr} &\EV^{\hat y}(T_{\DRc}) \notag\\
 & \leq \left(1 + O\left(\frac{r}{R}\right)\right) \min_{\hat y \in \Dr} \EV^{\hat y}(T_{\DRc}). \label{eq:BoundEscapeTimeToRStar}
\end{align}
For the sum in \eqref{eq:ExpValTau}, the Harnack inequality 
\eqref{eq:IntHarn2Toral} yields, for any $\hat y, \hat y' \in \Dr$, 
\begin{align}
 \sum_{\hat z \in \DRc} & \hat{H}_{\DRc}(\hat y,\hat z) \EV^{\hat z}(T_{\Dr}) \label{eq:UseHarnack}\\
 & = \left(1 + O\left(\frac{r}{R}\right)\right)\sum_{\hat z \in \DRc} \hat{H}_{\DRc}(\hat y',\hat z) \, \EV^{\hat z}(T_{\Dr}) \notag\\
 &  + O(R^{-M} \log R \lor K^{-M} R^2) \sum_{\hat z \in \hp(D(x,R+s)^c_K)} \EV^{\hat z}(T_{\Dr}). \notag
\end{align} 
The last term of \eqref{eq:UseHarnack} is zero if $p_1$ is finite range, by taking $s$ large enough so, due to \eqref{eq:IntHarn2Toral}, the error term does not appear. Otherwise, the sum needs to be controlled: since $R = cK^{1-\epsilon}$ and $\epsilon \geq 0$ is small, the Harnack inequality error is bounded above by 
\begin{align*}
cR^{-M} \log R & = c'K^{-4-2\beta+\epsilon(4+2\beta)} \log K \ll cK^{-M} R^2 & = cK^{-4-2\beta+2-2\epsilon} = cK^{-2-2\beta-2\epsilon} 
\end{align*}
and by 
\eqref{eq:ExTbdDnsBounds} with $R = cK^{1-\epsilon}$, the sum is bounded by $cK^{4-2\epsilon} (\log K)^2$.
Together these, with \eqref{eq:DRcBounds} and \eqref{eq:BoundEscapeTimeToRStar}, bound the last term of \eqref{eq:UseHarnack}: 
\begin{align}
c(R^{-M} &\log R \lor K^{-M} R^2)\sum_{\hat z \in \DRc} \EV^{\hat z}(T_{\Dr}) \label{eq:BoundExpAnnSum} \\
 & \leq cK^{2-2\beta-4\epsilon}(\log K)^2 \leq cK^{-2\beta-2\epsilon}(\log K)^2 \min_{\hat y \in \Dr} \EV^{\hat y}(T_{\DRc}). \notag
\end{align} 
Hence, by \eqref{eq:ExpValTau}-\eqref{eq:BoundExpAnnSum}, 
\begin{align}  
\max_{\hat y \in \Dr} \EV^{\hat y}(\tau)
\leq \bigg(1 & + O\left(\frac{r}{R}\right) + O(s^{-1}) \label{eq:TauHarnack}\\
 & + O\left(K^{-2\beta-2\epsilon}(\log K)^2\right)\bigg) \min_{\hat y \in \Dr} \EV^{\hat y}(\tau). \notag
\end{align}
Taking also $c_1 \geq 3c_0$, we get \eqref{eq:ExcursionLemma1} by combining \eqref{eq:OuterExcursionBounds} and \eqref{eq:TauHarnack}. $\qed$

The next corollary gives upper bounds for the hitting time of $\Dr$, and improves on \eqref{eq:ToralDiscEntryExpectedTime} for certain large radii.

\begin{cor} \label{lem:LPLemma3.1eqns23}
With the same hypotheses as above, 
\begin{align} 
\max_{\hat x \in \Z^2_K} \, \max_{\hat w \in \hp(\bd D(x,R)_R)} \EV^{\hat w}(T_{\Dr}) \leq c_1 K^2 \log\left(\frac{R}{r}\right); \label{eq:ExcursionLemma2} \\
\max_{\hat x \in \Z^2_K} ||T_{\Dr}|| \leq c_1 K^2 \log\left(\frac{K}{r}\right). \label{eq:ExcursionLemma3}
\end{align}
\end{cor}

\pf Consider \eqref{eq:ExpValTau} for $\hat y \in \Dr$ escaping to $\hp(D(x,4R)^c_K)$  instead of $\DRc$, before returning. Then, by \eqref{eq:ExcursionLemma1}, 
\begin{align} \label{eq:ExpValTau4m}
\sum_{\hat z \in \hp(D(x,4R)^c_K)} & \hat{H}_{\hp(D(x,4R)^c_K)}(\hat y,\hat z) \, \EV^{\hat z}(T_{\Dr}) \notag \\
 & \leq cK^2 \log(4R/r) \leq c'K^2 \log(R/r).
\end{align}
Using the strong Markov property at $T_{\hp(D(x,4R)^c_K)}$, \eqref{eqn:EscapeDiscExpTorus}, 
\eqref{eq:IntHarn2Toral}, \eqref{eq:ExpValTau4m}, 
\eqref{eq:ExTbdDnsBounds}, 
and \eqref{eq:BoundExpAnnSum}, 
we have for any $\hat w \in \hp(\bd D(x,R)_R)$ and some universal $c < \infty$, 
\begin{align} 
\EV^{\hat w}&(T_{\Dr}) \leq \EV^{\hat w}(T_{\hp(D(x,4R)^c_K)}) \label{eq:LP3.17ish}\\
 & \,\,\,\, + \EV^{\hat w}(T_{\Dr} - T_{\hp(D(x,4R)^c_K)}; T_{\Dr} > T_{\hp(D(x,4R)^c_K)}) \notag\\
 & \leq c\left[ (4R+1)^2 + \sum_{\hat z \in \hp(D(x,4R)^c_K)} \hat{H}_{\hp(D(x,4R)^c_K)}(\hat w,\hat z) \, \EV^{\hat z}(T_{\Dr})\right] \notag
\end{align}
\begin{align*} 
 & \leq c\bigg[ (4R+1)^2 + \sum_{\hat z \in \hp(D(x,4R)^c_K)} \left[ \left(1 + O\left(\frac{r}{R}\right)\right)\hat{H}_{\hp(D(x,4R)^c_K)}(\hat y,\hat z) \right] \EV^{\hat z}(T_{\Dr}) \notag\\
 & \,\,\,  + O(R^{-M} \log R \lor K^{-M} R^2) \sum_{\hat z \in \hp(D(x,4R+s)^c_K)} \EV^{\hat z}(T_{\Dr}) \bigg] 
 \leq c K^2 \log(R/r). \notag
\end{align*}
Setting $c_1 \geq c$, we have \eqref{eq:ExcursionLemma2}.  \eqref{eq:ExcursionLemma3} follows directly from \eqref{eq:ExcursionLemma2}, 
by considering $S$ projected onto $\Z^2_{24K}$ instead of $\Z^2_K$ for the furthest-out points $\hat w$. Note that, for these $\hat w$ such that $|\hat w - \hat x| > \frac{K}{24}$ on $\Z^2_K$, \eqref{eq:ExcursionLemma2} on $\Z^2_{24K}$ and the fact that annulus entrance takes longer on larger spaces, 
\[ \EV^{\hat w}(T_{\Dr}) \leq \EV^{\hat w}(T_{\hat \pi_{24K}(\bd D(x,r)_s)}) \leq c (24K)^2 \log(24K/r) \leq c_1 K^2 \log(K/r). \qed \]

\subsection{Decoupling an excursion from its endpoints} \label{sec:xband} 

\newcommand{\hDxrnksk}{\hp(\bd D(x,r_{n,k})_{s_k})}
\newcommand{\hDxrnkP}{\hp(D(x,r_{n,k}'))}
\newcommand{\hDxrnk}{\hp(D(x,r_{n,k}))}
\newcommand{\hDrnlpslp}{\hp(\bd D(0,r_{n,l+1})_{s_{l+1}})}
\newcommand{\hDrnlsl}{\hp(\bd D(0,r_{n,l})_{s_l})}
\newcommand{\hDrnlmslm}{\hp(\bd D(0,r_{n,l-1})_{s_{l-1}})}
\newcommand{\hDrnlc}{\hp(D(0, r_{n,l})^c_K)}
\newcommand{\hDrnlpc}{\hp(D(0, r_{n,l+1})^c_K)}

Let $n>13$ and set the following variables as defined in \eqref{eq:rsKdefn}: \index{Kn@$K_n$}
\begin{align*}
r_{n,k} = e^n n^{3k}, s_k = n^4, & & r_{n,k}' = r_{n,k}+s_{k},  & & k=0,1, \ldots, n; \notag\\
s_{n-1}^{n \downarrow} = \sqrt{r_{n,n-1}}
\end{align*}
and set $K_n := n^{\gbar} r_{n,n}$, where $\gbar \in [b,b+4]$ for some $b = b(p_1) \geq 10$, to be determined in Section \ref{ch:Late}. 

We say that, for a point $\hat x \in \Z^2_K$, and a path $\omega$ starting at $\hat x_0 \in \Z^2_{K_n}$, $\hat x_0 \neq \hat x$, the path $\omega$ \textbf{does not skip $\hat x$-bands} \index{does not skip $\hat x$-bands} if the path's  entrances and exits from the $r_{n,k}$-sized concentric discs around $\hat x$ are made by small or annulus-targeted jumps, not by medium or large untargeted jumps. 
More formally, a path does not skip $\hat x$-bands for a specified period of time if, during that time, escapes from $\hDxrnk$ 
 and entrances to $\hDxrnkP$ land in $\hDxrnksk$\footnote{ That is, with the exception of level $n-1$: entrances to $\hp(D(x,r_{n,n-1}+s_{n-1}^{n \downarrow}))$ land in the thicker band $\hp(\bd D(x,r_{n,n-1})_{s_{n-1}^{n \downarrow}})$. This is for the purposes of re-entering the level structure from the outermost level $n$; see \eqref{eq:bn0} for details, and assume this notation for excursions from level $n$ down to level $n-1$ if it is not mentioned.}.

By the strong Markov property, the only effect that one excursion between annuli has on another is via its beginning and ending points. In this section we build a structure in which to analyze the dependence on these endpoints for a special class of excursions.

The excursions we wish to examine are those from inside $\hp(D(0, r'_{n,l-1}))$ out to $\hDrnlc$ prior to ``one larger'' disc escape at $T_{\hDrnlpc}$. Consider a random path starting between these sets at $\hat z \in \hp(\bd D(0, r_{n,l})_{s_l})$. Focusing on annulus-based excursion end points $\hat w \in \hDrnlpslp$ and $l$ large, let $\mathcal{H}_{n,l-1 \uparrow l}$ \index{Hnl@$\mathcal{H}_{n,l-1  \uparrow l}$} be the $\sigma$-algebra of outward excursions $\hp(D(0, r_{n,l-1}')) \to \hDrnlc$ prior to $T_{\hDrnlpc}$. Let $\tau_0 = 0$, and for 
$i=0,1,2,\ldots$, define the excursion endpoint times 
\[\begin{array}{lll}
\tau_{2i+1} & = & \inf\{ k \geq \tau_{2i} : \hat{S}_k \in \hp(D(0, r_{n,l-1}')) \cup \hDrnlpc\}\\
\tau_{2i+2} & = & \inf\{ k \geq \tau_{2i+1} : \hat{S}_k \in \hDrnlc\}.\\
\end{array}\]
Abbreviating $\overline{\tau} = T_{\hDrnlpc}$, note that $\overline{\tau} = \tau_{2I+1}$ for some (unique) non-negative integer $I$. Then $\mathcal{H}_{n,l-1 \uparrow l}$ is the $\sigma$-algebra generated by the excursions $\{\hat{e}_{(j)}: j=1,\ldots,I\}$, where $\hat{e}_{(j)} = \{\hat{S}_k: \tau_{2j-1} \leq k \leq \tau_{2j}\}$ is the $j$th excursion $\hp(D(0, r_{n,l-1}')) \to \hDrnlc$. (The event $\{I=0\}$ is, of course, also included.)

\begin{figure}[!ht]
  \centering
    \includegraphics[width=4in]{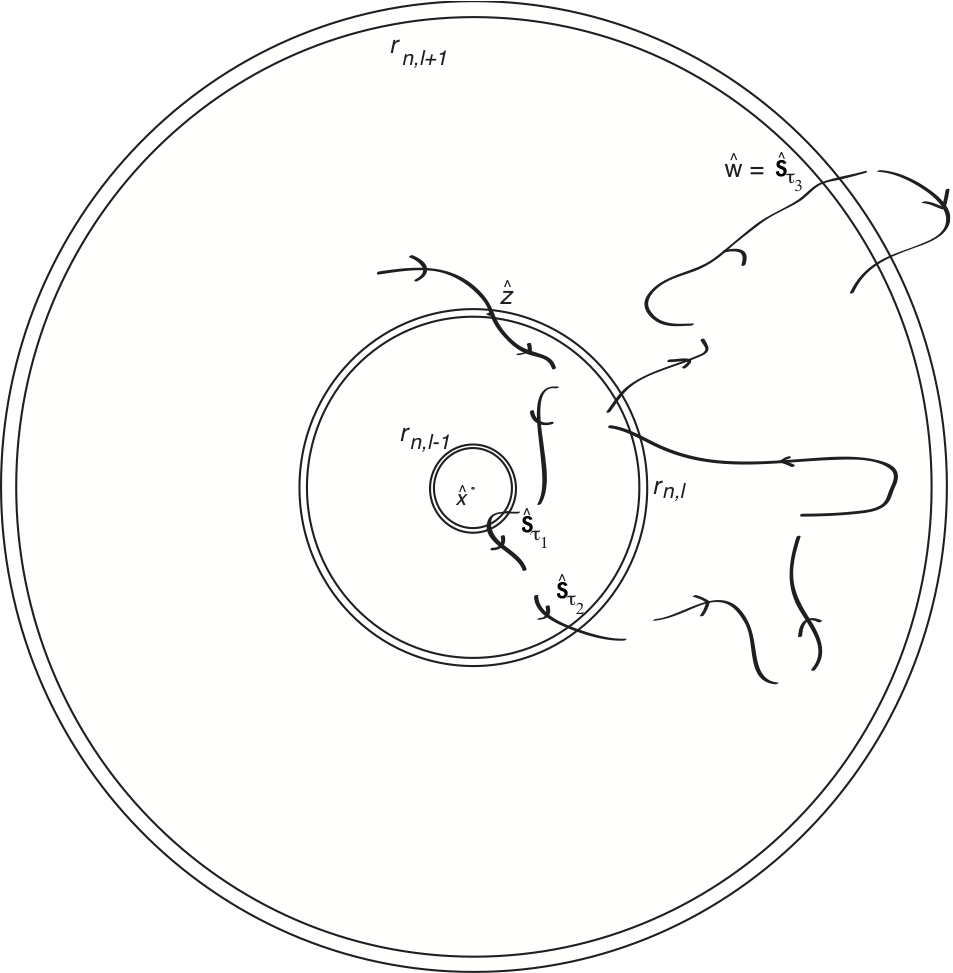}
    \parbox{4in}{
    \caption{Sample excursions - $\hat{e}_{(1)}$ is between $\hat{S}_{\tau_1}$ and $\hat{S}_{\tau_2}$. $I = 2$ for this path.}
    \label{fig:three_annuli_and_excursions_between}}
\end{figure}

Let $\mathcal{F}_{j} = \sigma(\hat{S}_k: k=0,1,\ldots,j)$, and for any stopping time $\tau$, let $\mathcal{F}_{\tau}$ denote the collection of events $A$ \st $A \cap \{\tau = j\} \in \mathcal{F}_j$ for all $j$.

We will focus on paths which do not skip $\hat x$-bands over a number of concentric annulus excursions. Let  $\Omega_{\hat x,n,l+1,m}^{i-1,\ldots,j}$ \index{Oxnim@$\Omega_{\hat x,n,i,m}^{i-1,\ldots,j}$} denote the set of paths which do not skip $\hat x$-bands on excursions between levels $k=i-1, i, \ldots, j$ until completion of the first $m$ outward excursions from $\hp(D(x, r_{n,l}')) \to \hp(D(x, r_{n,l+1})^c_K)$, and $\Omega_{\hat x,n,l+1,m}^{A}$ the same for the levels in the index set $A$.
Our first lemma shows that excursion paths faithful to hitting $\hat x$-bands are ``almost'' independent of their beginning and ending points.

\begin{lem} \label{lem:FPLemma8.1} 
Uniformly in $l$, $n$, $K_n$, $B_n \in \mc{H}_{n,l-1 \uparrow l}$, $\hat z, \hat z' \in \hDrnlsl$, and $\hat w \in \hDrnlpslp$, 
\begin{align} 
P^{\hat z}\bigg(B_n \cap \Omega_{\hat 0,n,l+1,1}^{l-1,l,l+1}  & \, \bigg| \, \hat{S}_{T_{\hp(D(0,r_{n,l+1})^c_K)}} = \hat w\bigg) \label{eq:FPLemma8.1(8.1)} \\ 
 & = (1 + O(n^{-3})) P^{\hat z} \left(B_n \cap \Omega_{\hat 0,n,l+1,1}^{l-1,l,l+1} \right) \notag
\end{align}
and 
\begin{equation} \label{eq:FPLemma8.1(8.2)} 
P^{\hat z} \left(B_n \cap \Omega_{\hat 0,n,l+1,1}^{l-1,l,l+1} \right) = (1 + O(n^{-3} \log n)) P^{\hat z'} \left(B_n \cap \Omega_{\hat 0,n,l+1,1}^{l-1,l,l+1} \right).
\end{equation}
\end{lem}

\pf Fixing a starting point $\hat z \in \hDrnlsl$, it suffices to consider $B_n \in \mc{H}_{n,l-1 \uparrow l}$ \st $P^{\hat z}(B_n) > 0$. Fix such a set $B_n$ and an ending point $\hat w \in \hDrnlpslp$. Using the notation just introduced, for any $i \geq 1$, we can write 
\[\begin{array}{l}
B_n \cap \Omega_{\hat 0,n,l+1,1}^{l-1,l,l+1} \cap \{I = i\}\\
\,\,\,\, = B_{n,i} \cap A_i \cap \{\tau_{2i} < \overline{\tau}\} 
  \cap (\{I = 0, \hat{S}_{\overline{\tau}} \in \hDrnlpslp\} \circ \theta_{\tau_{2i}})
\end{array}\]
for some $B_{n,i} \in \mc{F}_{\tau_{2i}}$, where 
\[A_i = \{ \hat{S}_{\tau_{2j-1}} \in \hDrnlmslm, \, \hat{S}_{\tau_{2j}} \in \hp(\bd D(0,r_{n,l})_{s_l})\, ,\forall j \leq i \} \in \mc{F}_{\tau_{2i}},\]
so by the strong Markov property at $\tau_{2i}$, 
\begin{align*}
& P^{\hat z}( \{\hat{S}_{\overline{\tau}} = \hat w\} \cap B_n \cap \Omega_{\hat 0,n,l+1,1}^{l-1,l,l+1} \cap \{I = i\}) \notag\\
& \,\,\,\, = \EV^{\hat z}[ P^{\hat{S}_{\tau_{2i}}}(\hat{S}_{\overline{\tau}} = \hat w; \, I = 0); B_{n,i} \cap A_i \cap \{\tau_{2i} < \overline{\tau}\} ]; \notag\\
& P^{\hat z}(B_n \cap \Omega_{\hat 0,n,l+1,1}^{l-1,l,l+1} \cap \{I = i\}) \notag\\
& \,\,\,\, = \EV^{\hat z}[ P^{\hat{S}_{\tau_{2i}}}(\hat{S}_{\overline{\tau}} \in \hDrnlpslp; \, I = 0); 
  B_{n,i} \cap A_i \cap \{\tau_{2i} < \overline{\tau}\} ]. \notag
\end{align*}
Consequently, for all $i \geq 1$, 
\begin{align}
& P^{\hat z}( \{\hat{S}_{\overline{\tau}} = \hat w\} \cap B_n \cap \Omega_{\hat 0,n,l+1,1}^{l-1,l,l+1} \cap \{I = i\} ) \label{eq:FP(8.3)} \\
 & \,\,\,\, \geq P^{\hat z}( B_n \cap \Omega_{\hat 0,n,l+1,1}^{l-1,l,l+1} \cap \{I = i\} ) 
   \min_{\hat x \in \hp(\bd D(0,r_{n,l})_{s_l})} \frac{P^{\hat x}(\hat{S}_{\overline{\tau}} = \hat w; I = 0)}{P^{\hat x}(\hat{S}_{\overline{\tau}} \in \hDrnlpslp; I = 0)}. \notag
\end{align}

Note that 
\[\{I = 0\} = \{\tau = T_{\hDrnlpc} < T_{\hp(D(0, r_{n,l-1}'))}\}.\]
Necessarily, $P^{\hat z}(B_n | I=0) \in \{0,1\}$ and is independent of $\hat z$ for any $B_n \in \mc{H}_{n,l-1 \uparrow l}$, implying that (\ref{eq:FP(8.3)}) applies for $i=0$ as well. Hence, by \eqref{eq:InteriorHarnackZoomIn} and  \eqref{eq:InteriorHarnack}, 
there exists $c < \infty$ \st for any $\hat z, \hat x \in \hp(\bd D(0,r_{n,l})_{s_l})$ and $\hat w \in \hDrnlpslp$,
\[\frac{P^{\hat x}(\hat{S}_{\overline{\tau}} = \hat w; I = 0)}{P^{\hat x}(\hat{S}_{\overline{\tau}} \in  \hDrnlpslp; I = 0)}
 \geq (1 - cn^{-3}) \hat{H}_{\hDrnlpc}(\hat z,\hat w).\]
We note that, since \eqref{eq:InteriorHarnackZoomIn} and \eqref{eq:InteriorHarnack} accommodate starting points up to a square root of 
 the distance away from their level's starting radius of $r_{n,l}$, this bound is good for even the wide band $s_{n-1}^{n \downarrow} = \sqrt{r_{n,n-1}} \ll r_{n,n-1}$ as a starting point (this is the case $l = n-1$). 

  
Hence, summing \eqref{eq:FP(8.3)} over $I=0,1,\ldots$, we get that 
\begin{equation*} 
P^{\hat z}(\{\hat{S}_{\overline{\tau}} = \hat w\} \cap B_n \cap \Omega_{\hat 0,n,l+1,1}^{l-1,l,l+1}) 
   \geq (1 - cn^{-3}) P^{\hat z}(B_n \cap \Omega_{\hat 0,n,l+1,1}^{l-1,l,l+1}) \hat{H}_{\hDrnlpc}(\hat z,\hat w).
\end{equation*}
A similar argument shows that 
\begin{equation*} 
P^{\hat z}(\{\hat{S}_{\overline{\tau}} = \hat w\} \cap B_n \cap \Omega_{\hat 0,n,l+1,1}^{l-1,l,l+1}) 
   \leq (1 + cn^{-3}) P^{\hat z}(B_n \cap \Omega_{\hat 0,n,l+1,1}^{l-1,l,l+1}) \hat{H}_{\hDrnlpc}(\hat z,\hat w),
\end{equation*}
and we obtain \eqref{eq:FPLemma8.1(8.1)}.

By the strong Markov property at $\tau_1$, for any $\hat z \in \hp(\bd D(0,r_{n,l})_{s_l})$, 
\begin{align*}
& P^{\hat z}(B_n \cap \Omega_{\hat 0,n,l+1,1}^{l-1,l,l+1}) = P^{\hat z}(B_n \cap \Omega_{\hat 0,n,l+1,1}^{l-1,l,l+1} \cap \{I=0\}) \notag\\
& \,\,\,\, + \sum_{\hat x \in \hDrnlmslm} \hat{H}_{\hp(D(0,r_{n,l-1}')) \cup \hp(D(0,r_{n,l+1})^c)}(\hat z,\hat x) P^{\hat x}(B_n \cap \Omega_{\hat 0,n,l+1,1}^{l-1,l,l+1}). \notag
\end{align*}
The first term is handled by (\ref{eq:ToralGamblersSuccess}). \eqref{eq:FPLemma8.1(8.2)} then follows from \eqref{eq:ExteriorHarnackToralZoomOut}. $\qed$ 

Next, we examine excursions going inward: let $\mathcal{G}_{n,l+1 \downarrow l}^{\hat x}$ \index{Gnlx@$\mathcal{G}_{n,l+1 \downarrow l}^{\hat x}$} denote the $\sigma$-algebra of excursions from $\hp(D(x, r_{n,l+1})^c_{K_n})$ into $\hp(D(x, r_{n,l}'))$. 
To this end, let $\hat x \in \Z^2_{K_n}$, let $\overline{\tau}_0 = 0$ and for $i=1,2,\ldots$ define 
\[\begin{array}{lll}
\tau_i & = & \inf\{ k \geq \overline{\tau}_{i-1}: \hat{S}_k \in \hp(D(x, r_{n,l}'))\},\\
\overline{\tau}_i & = & \inf\{ k \geq \tau_i: \hat{S}_k \in \hp(D(x, r_{n,l+1})^c_{K_n})\}.
\end{array}\]
Then $\mathcal{G}_{n,l+1 \downarrow l}^{\hat x}$ is the $\sigma$-algebra generated by the excursions $\{\hat{e}^{(j)}: j=1,\ldots\}$, where $\hat{e}^{(j)} = \{\hat{S}_k: \overline{\tau}_{j-1} \leq k \leq \tau_{j}\}$ is the $j$th excursion  $\hp(D(x, r_{n,l+1})^c_{K_n}) \to \hp(D(x, r_{n,l}'))$ (so for $j=1$ we begin at $t=0$).

Let $\mc{H}_{n,l-1 \uparrow l}^{\hat x}(m)$ \index{Hnlxm@$\mc{H}_{n,l-1 \uparrow l}^{\hat x}(m)$} be the $\sigma$-algebra of excursions from $\hp(D(x, r_{n,l-1}'))$ out to $\hp(D(x, r_{n,l})^c_{K_n})$ during the first $m$ excursions from $\hp(D(x, r_{n,l}'))$ out to $\hp(D(x, r_{n,l+1})^c_{K_n})$, \ie from $\tau_1$ to $\overline{\tau}_m$. In more detail, for each $j=1,2,\ldots,m$, let $\overline{\zeta_{j,0}} = \tau_j$ and for $i=1,\ldots$, define 
\[\begin{array}{lll}
\zeta_{j,i} & = & \inf\{ k \geq \overline{\zeta}_{j,i}: \hat{S}_k \in \hp(D(x, r_{n,l-1}'))\},\\
\overline{\zeta}_{j,i} & = & \inf\{ k \geq \zeta_{j,i}: \hat{S}_k \in \hp(D(x, r_{n,l})^c_{K_n})\},\\
v_{j,i} & = & \{\hat{S}_k: \zeta_{j,i} \leq k \leq \overline{\zeta_{j,i}}\},\\
Z^j & = & \sup \{i \geq 0: \overline{\zeta}_{j,i} < \overline{\tau}_j\}.
\end{array}\]
Then $\mc{H}_{n,l-1 \uparrow l}^{\hat x}(m)$ is the $\sigma$-algebra generated by the intersection of the $\sigma$-algebras $\mc{H}_{n,l,j}^{\hat x} = \sigma(v_{j,i}: i=1,2,\ldots,Z^j)$ \index{Hnljx@$\mc{H}_{n,l,j}^{\hat x}$} of the excursions between $\tau_j$ and $\overline{\tau}_j$, for $j=1,2,\ldots,m$.

\begin{lem} \label{lem:FPLemma8.2} 
There exists $C < \infty$ such that, uniformly over all $m \leq (n \log n)^2$, $l, \hat x \in \Z^2_{K_n}$ and $\hat y_0, \hat y_1 \in \Z^2_{K_n} \setminus \hp(D(x,r_{n,l}'))$, and $H \in \mc{H}_{n,l-1 \uparrow l}^{\hat x}(m)$,
\begin{equation} \label{eq:FP(8.6)}
\begin{array}{l}
(1 - Cmn^{-3} \log n) P^{\hat y_1}(H \cap \Omega_{\hat x,n,l+1,m}^{l-1,l,l+1})\\
\,\,\,\, \leq P^{\hat y_0}(H \cap \Omega_{\hat x,n,l+1,m}^{l-1,l,l+1} | \mathcal{G}_{n,l+1 \downarrow l}^{\hat x})
 \leq (1 + Cmn^{-3} \log n) P^{\hat y_1}(H \cap \Omega_{\hat x,n,l+1,m}^{l-1,l,l+1}).
\end{array}
\end{equation}
\end{lem}

\pf Applying the Monotone Class Theorem to the algebra of their finite disjoint unions, it suffices to prove (\ref{eq:FP(8.6)}) for the generators of the $\sigma$-algebra $\mc{H}_{n,l-1 \uparrow l}^{\hat x}(m)$ of the form 
$H = H_1 \cap H_2 \cap \cdots \cap H_m$, 
 with $H_j \in \mc{H}_{n,l,j}^{\hat x}$ for $j=1,\ldots,m$. Conditioned upon $\mc{G}_{n,l+1 \downarrow l}^{\hat x}$, the events $H_j$ are independent. Further, each $H_j$ then has the conditional law of an event $B_j$ in the $\sigma$-algebra $\mc{H}_{n,l-1 \uparrow l}$ of Lemma \ref{lem:FPLemma8.1}, for some random end points $\hat z_j = \hat{S}_{\tau_j} - \hat x \in  \hDrnlsl$ and $\hat w_j = \hat{S}_{\overline{\tau}_j} - \hat x \in \hDrnlpslp$, both measurable on $\mc{G}_{n,l+1 \downarrow l}^{\hat x}$. By our conditions, the uniform estimates (\ref{eq:FPLemma8.1(8.1)}) and (\ref{eq:FPLemma8.1(8.2)}) yield that for any fixed $\hat z' \in \hDrnlsl$,
\[\begin{array}{lll}
P^{\hat y_0}(H \cap \Omega_{\hat x,n,l+1,m}^{l-1,l,l+1} | \mc{G}_{n,l+1 \downarrow l}^{\hat x})
 & = & P^{\hat y_0}(\cap_{j=1}^m (H_j \cap \Omega_{\hat x,n,l+1,1}^{l-1,l,l+1}) | \mc{G}_{n,l+1 \downarrow l}^{\hat x})\\
 & = & \prod_{j=1}^m P^{\hat z_j}(B_j \cap \Omega_{\hat x,n,l+1,1}^{l-1,l,l+1} | \hat{S}_{T_{D(0,r_{n,l})^c}} = \hat w_j)\\
 & = & \prod_{j=1}^m (1 + O(n^{-3})) P^{\hat z_j}(B_j \cap \Omega_{\hat x,n,l+1,1}^{l-1,l,l+1})\\
 & = & (1 + O(n^{-3} \log n))^m \prod_{j=1}^m P^{\hat z'}(B_j \cap \Omega_{\hat x,n,l+1,1}^{l-1,l,l+1}).
\end{array}\]
Since $m \leq (n \log n)^2$ and the last expression above neither depends on $\hat y_0 \in \Z^2_K$ nor on the extra information in $\mc{G}_{n,l+1 \downarrow l}^{\hat x}$, we get \eqref{eq:FP(8.6)}. $\qed$

Now that we have control over the excursion structure of paths that do not skip $\hat x$-bands, we will control their layered excursion counts. Fix $0 < a < 2$, and define \index{vn@$v_n(a)$} $v_k = v_k(a) := 3ak^2 \log k$ for $k=2,3,...,n$, and $N^{\hat x}_{n,l}$, $l=2,\ldots,n-1$, \index{Nxnl@$N^{\hat x}_{n,l}$} as the number of excursions from $\hp(D(x,r_{n,l-1}'))$ out to $\hp(D(x,r_{n,l})^c_{K_n})$ until time \index{Rxn@$\mc{R}_n^{\hat x}$} $\mc{R}_n^{\hat x}(a)$, the time that $v_n$ excursions from  $\hp(D(x,r_{n,n-1}))$ out to $\hp(D(x,r_{n,n})^c_{K_n})$ have been completed. Let $m \stackrel{k}{\sim} v$ \index{simk@$\stackrel{k}{\sim}$} denote the bound $|m-v| \leq k$. Finally, let $N^{\hat x}_{n,0}$ be the number of visits to $\hat x$ before $\mc{R}_n^{\hat x}(a)$.

\begin{lem} \label{lem:FPCor8.3} 
Let \index{Gnly@$\Gamma_{n,l}^{\hat y}$} $\Gamma_{n,l}^{\hat y} := \{N_{n,i}^{\hat y} = m_i: i=0,2,\ldots,l-1\} \cap \Omega_{\hat y,n,l+1,m_l}^{1,\ldots,l}$. Then, for any $1 < n_0 < n$, uniformly over all $n_0 \leq l \leq n-1$, $m_l \stackrel{l}{\sim} v_l$, $\{m_i: i=0,2,\ldots,l\}$, $\hat y \in \Z^2_{K_n}$, and $\hat x_0, \hat x_1 \in \Z^2_{K_n} \setminus \hp(D(y,r_{n,l}'))$, 
\begin{align} 
 & P^{\hat x_0}(\Gamma_{n,l}^{\hat y}, N_{n,l}^{\hat y} = m_l | \mc{G}_{n,l \downarrow l-1}^{\hat y}) \notag\\
 & \,\,\,\, = ( 1 + O(n^{-1} (\log n)^2) ) P^{\hat x_1}(\Gamma_{n,l}^{\hat y} | N_{n,l}^{\hat y} = m_l) 1_{\{N_{n,l}^{\hat y} = m_l\}}. \label{eq:FP(8.8)}
\end{align}
\end{lem}

\pf For $j=1,2,\ldots$ and $i=2,\ldots,l$, let $Z^j_i$ \index{Zji@$Z^j_i$} denote the number of excursions from $\hp(D(x,r_{n,i}'))$ out to $\hp(D(x,r_{n,i+1})^c_{K_n})$ by the random walk during the time interval $[\tau_j,\overline{\tau_j}]$. The event 
\[H = \left\{ \sum_{j=1}^{m_l} Z_i^j = m_i : i=2,\ldots,l-1 \right\} \cap \Omega_{\hat y,n,l,m_{l-1}}^{2,\ldots,l-1}\]
belongs to the $\sigma$-algebra $\mc{H}_{n,l-1 \uparrow l}^{\hat y}(m_l)$ of Lemma \ref{lem:FPLemma8.2}. It is easy to verify that, starting at any $\hat x_0 \not \in \hp(D(y, r_{n,l}'))$, when the event $\{N_{n,l}^{\hat y} = m_l\} \in \mc{G}_{n,l \downarrow l-1}^{\hat y}$ occurs, it implies that $N_{n,i}^{\hat y} = \sum_{j=1}^{m_l} Z_i^j$ for $i=2,\ldots,l$. Thus, setting $H' = H \cap \Omega_{\hat y,n,l+1,m_l}^{l-1,l,l+1}$, 
\begin{equation} \label{eq:FP(8.9)}
P^{\hat x_0}(\Gamma_{n,l}^{\hat y} | \mc{G}_{n,l \downarrow l-1}^{\hat y}) 1_{\{N_{n,l}^{\hat y} = m_l\}}
 = P^{\hat x_0}(H'| \mc{G}_{n,l \downarrow l-1}^{\hat y}) 1_{\{N_{n,l}^{\hat y} = m_l\}}.
\end{equation}
With $m_l/(l^2 \log l)$ bounded above, by (\ref{eq:FP(8.6)}) we have, uniformly in $\hat y \in \Z^2_{K_n}$ and $\hat x_0, \hat x_1 \in \Z^2_{K_n} \setminus \hp(D(y,r_{n,l}'))$, 
\begin{equation} \label{eq:FP(8.10)}
P^{\hat x_0}(H' | \mc{G}_{n,l \downarrow l-1}^{\hat y}) = (1 + O(n^{-1} (\log n)^2)) P^{\hat x_0}(H').
\end{equation}
Hence, 
\begin{equation} \label{eq:FP(8.11)}
\begin{array}{l}
P^{\hat x_0}(\Gamma_{n,l}^{\hat y} | \mc{G}_{n,l \downarrow l-1}^{\hat y}) 1_{\{N_{n,l}^{\hat y} = m_l\}} = (1 + O(n^{-1} (\log n)^2)) P^{\hat x_1}(H') 1_{\{N_{n,l}^{\hat y} = m_l\}}.
\end{array}
\end{equation}
Setting $\hat x_0 = \hat x_1$ and taking expectations with respect to $P^{\hat x_0}$ yields 
\begin{align} 
P^{\hat x_1}(\Gamma_{n,l}^{\hat y} | N_{n,l}^{\hat y} = m_l)
 & = (1 + O(n^{-1} (\log n)^2)) P^{\hat x_1}(H') \label{eq:FP(8.13)}\\ 
\implies P^{\hat x_1}(\Gamma_{n,l}^{\hat y} | N_{n,l}^{\hat y} = m_l) 1_{\{N_{n,l}^{\hat y} = m_l\}}
 & = (1 + O(n^{-1} (\log n)^2)) P^{\hat x_1}(H') 1_{\{N_{n,l}^{\hat y} = m_l\}}  \notag\\
 & = (1 + O(n^{-1} (\log n)^2)) P^{\hat x_0}(\Gamma_{n,l}^{\hat y} | \mc{G}_{n,l \downarrow l-1}^{\hat y}) 1_{\{N_{n,l}^{\hat y} = m_l\}} \notag 
\end{align}
where we used (\ref{eq:FP(8.11)}) for the last equality. With $\{N_{n,l}^{\hat y} = m_l\} \in \mc{G}^y_{n,l \downarrow l-1}$, this is \eqref{eq:FP(8.8)}. $\qed$


\section{Late Points} \label{ch:Late}

We define the \emph{cover time} \index{cover time} of $\Z^2_K$ by the random walk $\hat{S}$ to be the maximum first visiting time over all points in $\ZK$: if $\mc{T}_K(\hat x) = \inf\{t \geq 0: \hat{S}_t = \hat x\}$ is the first time visiting $\hat x$, then the cover time of $\ZK$ is
\begin{equation} \label{eqn:CoveringTime}
\mc{T}_{cov}(\Z^2_K) := \max_{\hat x \in \ZK} \, \mc{T}_K(\hat x).
\end{equation}

In \cite{DPRZ2004}, Dembo, Peres, Rosen, and Zeitouni showed that the cover time of $\Z_K^2$ for simple random walk is asymptotic to ${4 \over \pi}(K \log K)^2$ as $K \to \infty$. This result was found via strong approximation techniques to Brownian motion. The team reproduced this result via purely random walk methods in \cite{DPRZ2006}, along with a multifractal analysis of the late points of the torus. Here we generalize results from \cite{BRFreq} and \cite{DPRZ2006} to gain similar results for toral random walks with jumps of infinite range.

Let $\alpha \in (0,1)$. Anticipating the result, we call $\hat x$ an $\alpha,K$-\emph{late point} \index{late point} of the random walk $\hat{S}$ on $\Z^2_K$ if $\mc{T}_K(\hat x) \geq \frac{4\alpha}{\pi_{\Gamma}}(K \log K)^2$. Set $\mc{L}_K(\alpha)$ \index{LK@$\mc{L}_K(\alpha)$} to be the set of $\alpha,K$-late points in $\Z^2_K$, \ie 
\[\mc{L}_K(\alpha) := \left\{ \hat x \in \Z^2_K: \frac{\mc{T}_K(\hat x)}{(K \log K)^2} \geq \frac{4\alpha}{\pi_{\Gamma}} \right\}. \]


\subsection{Upper bound of late point probabilities}

First we show that excursion times are concentrated around their mean, and relate excursions to hitting times.

\begin{lem} \label{lem:BoundSumOfExcLengths} 
 With the notation of Lemma \ref{lem:LPLemma3.1}, we can find $\de_0 > 0$ and $C>0$ such that, 
 if $R \leq K/24$ and $\de \leq \de_0$ with $\de \leq 6c_1(s^{-1} + r/R)$, then for all $\hat x, \hat x_0 \in \Z_K^2$,
\begin{equation} \label{eq:lemBoundExcTimesAbove}
 P^{\hat x_0}\left(\sum_{j=0}^{N} \tau^{(j)} \leq (1 - \de) N \frac{2K^2 \log(R/r)}{\pi_{\Gamma}}\right)
 \leq e^{-C\de^2 N (\log(R/r)/\log(K/r))}
\end{equation}
and
\begin{equation} \label{eq:lemBoundExcTimesBelow}
 P^{\hat x_0}\left(\sum_{j=0}^{N} \tau^{(j)} \geq (1 + \de) N \frac{2K^2 \log(R/r)}{\pi_{\Gamma}}\right)
 \leq e^{-C\de^2 N (\log(R/r)/\log(K/r))}.
\end{equation}
\end{lem}

\pf 
With $\tau=\tau^{(1)}=\lc T_{\DRc} + T_{\Dr}\circ\theta_{T_{\DRc}}\rc \circ\theta_{T_{\Dr}}$, 
\begin{align*}
 \max_{\hat y \in \Dr}  & \EV^{\hat y} (\tau^{n})
 \leq \max_{\hat y \in \Dr} \EV^{\hat y} \left(\lc T_{\DRc}
 + T_{\Dr}\circ\theta_{T_{\DRc}}\rc^{n}\right)\\
 & \leq \sum_{ j=0}^{ n} \binom{n}{j} \max_{y \in \Dr} 
     \EV^{\hat y} \big( T^{ j}_{\DRc}\,( \,T^{ n-j}_{\Dr}\circ\theta_{T_{\DRc}}) \big)\\
 & \leq \sum_{ j=0}^{ n} \binom{n}{j}
     \max_{\hat y \in \Dr} \EV^{\hat y} ( T^{ j}_{\DRc})
     \max_{\hat z \in \DRc} \EV^{\hat z} ( T^{ n-j}_{\Dr})\,.
\end{align*}
Let $u=\frac{2K^{ 2}}{\pi_{\Gamma}}\log(K/r)$ 
 and $u'=\frac{2K^{2}}{\pi_{\Gamma}}\log(R/r)$. 
Then, by (\ref{eq:KacMomentStrongMarkov}), (\ref{eq:ExcursionLemma2}), (\ref{eqn:EscapeDiscExpTorus}), and (\ref{eq:ExcursionLemma3}), we can bound the moments of $\tau$: there exist universal constants $c_1, c_2 < \infty$ \st for all $\hat x \in \Z^2_K$,
\begin{equation} \label{eq:BoundingTauMoments}
\begin{array}{lll}
\max_{\hat y \in \Dr}  \EV^{\hat y} (\tau^{n})
 & \leq & \max_{\hat y \in \Dr} \EV^{\hat y} (T_{\DRc})||T_{\DRc}||^{n-1} n!\\
 & & + 2 c_1 \sum_{j=0}^{n-1} n! ||T_{\DRc}||^j u' ||T_{\Dr}||^{n-j-1}\\
 & \leq & (n+1)! u' (c_2 u)^{n-1}.
\end{array}
\end{equation}
Taking $\eta = \delta/6 > 0$, with our choice of $r$ and $R$, it thus follows by (\ref{eq:ExcursionLemma1}) that for $\rho = c_3 u u'$ and all $\theta > 0$,
\begin{equation} \label{eq:TauLaplaceTransform}
\begin{array}{lll}
\max_{\hat x} \max_{\hat y \in \Dr}  \EV^{\hat y} (e^{-\theta \tau})
 & \leq & 1 - \theta \min_{\hat x} \min_{\hat y \in \Dr} \EV^{\hat y} (\tau)\\
 & & + \frac{\theta^2}{2} \max_{\hat x} \max_{\hat y \in \Dr} \EV^{\hat y} (\tau^2)\\
 & \leq & 1 - \theta(1-\eta)u' + \rho \theta^2\\
 & \leq & \exp(\rho \theta^2 - \theta(1-\eta)u').
\end{array}
\end{equation}
Since $\tau^{(0)} \geq 0$, using Markov's inequality, we bound the left-hand side of \eqref{eq:lemBoundExcTimesAbove} by 
\begin{align} \label{eq:ExcTimesAboveBound}
P^{\hat x_0}\left( \sum_{j=1}^N \tau^{(j)} \leq (1 - 6\eta) u'N \right)
 & \leq e^{\theta(1-3\eta)u'N} \EV^{\hat x_0}(e^{-\theta\sum_{j=1}^N \tau^{(j)}}) \\
 & \leq e^{-\theta u'N \delta/3} \left[ e^{\theta(1-\eta)u'} \max_{\hat y \in \Dr} \EV^{\hat y} (e^{-\theta \tau})  \right]^N, \notag
\end{align}
where the last inequality follows by the strong Markov property of $\hat{S}_t$ on $\{\mathfrak{T}_j\}$. Combining 
\eqref{eq:TauLaplaceTransform} and \eqref{eq:ExcTimesAboveBound} for $\theta = \delta u'/(6 \rho)$ results in 
\eqref{eq:lemBoundExcTimesAbove} for $C = 1/(36 c_3)$.

Since $\tau^{(0)} = T_{\Dr}$, by (\ref{eq:KacMomentStrongMarkov}) and (\ref{eq:ExcursionLemma3}), there exist universal constants $c_4$, $c_5 < \infty$ \st 
\[\max_{\hat x,\hat y} \EV^{\hat y}(e^{\tau^{(0)}/c_4 u}) \leq c_5.\]
This implies 
\[P^{\hat x_0}\left( \tau^{(0)} \geq \frac{\delta}{3} u'N \right)
  = P^{\hat x_0}\left( \frac{\tau^{(0)}}{c_4 u} \geq \frac{\delta}{3 c_4} \frac{u'}{u} N \right)
  \leq c_5 e^{(-3 c_4)^{-1} \delta (u'/u) N}.\]
Thus, the proof of \eqref{eq:lemBoundExcTimesBelow}, like in \eqref{eq:lemBoundExcTimesAbove}, comes down to bounding 
\[P^{\hat x_0}\left( \sum_{j=1}^N \tau^{(j)} \geq (1 + 4\eta) u'N \right)
 \leq e^{-\theta u'N \delta/3} \left[ e^{-\theta(1+2\eta)u'} \max_{\hat y \in \Dr} \EV^y (e^{\theta \tau}) \right]^N.\]
Noting that, by (\ref{eq:BoundingTauMoments}) and (\ref{eq:ExcursionLemma1}), there exists a universal constant $c_6 < \infty$ \st for $\rho = c_6 u u'$ and all $0 < \theta < 1/(2 c_2 u)$,
\begin{align}
\max_{\hat x} \max_{\hat y \in \Dr} \EV^{\hat y} (e^{\theta \tau})
 & \leq 1 + \theta \max_{\hat y \in \Dr} \EV^{\hat y} (\tau) + \sum_{n=2}^{\infty} \frac{\theta^n}{n!} \EV^{\hat y} (\tau^n) \notag\\
 & \leq 1 + \theta (1 + 2\eta) u' + \rho \theta^2 \\
 & \leq \exp(\theta (1 + 2\eta) u' + \rho \theta^2). \notag
\end{align}
Taking $\delta_0 < 3c_6 / c_2$, the proof of (\ref{eq:lemBoundExcTimesBelow}) now follows that of (\ref{eq:lemBoundExcTimesAbove}). $\qed$

Next we apply Lemma \ref{lem:BoundSumOfExcLengths} to bound the upper tail of $\mathcal{T}_K(\hat x)$, the first hitting time of $\hat x \in \Z^2_K$.

\begin{lem} \label{lem:HittingTimeTailProb}
 For any $\de > 0$ we can find $c < \infty$ and $K_0 < \infty$ such that,
 for all $K \geq K_0$, $b \geq 0$, and $\hat x, \hat x_0 \in \Z_K^2$,
\begin{equation} \label{eq:LemUpperBd}
 P^{\hat x_0}\left(\TT_K(\hat x) \geq b(K \log K)^2\right) \leq cK^{-(1-\de)\pi_{\Gamma} b/2}.
\end{equation}
\end{lem}
 
\pf Fix $\de \in (0,\de_0)$, where $\de_0$ is from Lemma \ref{lem:BoundSumOfExcLengths}. Let $R=\frac{K}{24}$ and $r=R/\log K$. Then Lemma \ref{lem:BoundSumOfExcLengths} applies for all $K \geq K_0$ and some $K_0 = K_0(\de) < \infty$. Fixing $b \geq 0$ and such $K$, let
\[n_K := (1 - \de)\frac{\pi_{\Gamma} b (\log K)^2}{2 \log (R/r)} = (1 - \de)\frac{\pi_{\Gamma} b (\log K)^2}{2 \log \log K}.\]
Then,
\begin{align}
P^{\hat x_0}\left(\TT_K(\hat x) \geq b(K \log K)^2\right)
 & \leq \,\, P^{\hat x_0}\left(\TT_K(\hat x) \geq \sum_{j=0}^{n_K} \tau^{(j)}\right) \notag\\
 & \,\, + P^{\hat x_0}\left(\sum_{j=0}^{n_K} \tau^{(j)} \geq b(K \log K)^2\right). \label{eq:SplitLateExcProbs}
\end{align}
The first probability in the sum in \eqref{eq:SplitLateExcProbs} is 
the probability of not hitting $\hat x$ during the first $n_K$ consecutive $\Dr \to \DRc \to \Dr$ excursions. By \eqref{eq:ProbZeroBeforeDisc}, 
\begin{align} \label{eq:HitDRcMissx}
P^{\hat x_1}\left(T_{\hat x} < T_{\DRc}\right)
 & = \left[ \frac{\log(R/r) + O(r^{-1/4})}{\log(R)} \right] \left( 1 + O(\log(R)^{-1}) \right)
\end{align}
uniformly for $\hat x_1 \in \Dr$. 
For any $\hat x_2 \in \DRc$, 
\begin{align} \label{eq:HitDrMissx}
P^{\hat x_2}\left(T_{\hat x} < T_{\Dr}\right) < 1. 
\end{align}
Hence, by \eqref{eq:HitDRcMissx} and \eqref{eq:HitDrMissx}, the first probability 
 in \eqref{eq:SplitLateExcProbs} is bounded above by 
\begin{align}
 & \max_{\stackrel{\hat x_1 \in \Dr}{\hat x_2 \in \DRc}} \bigg[\left( 1 - P^{\hat x_1}\left(T_{\hat x} < T_{\DRc}\right) \right) \left( 1 - P^{\hat x_2}\left(T_{\hat x} < T_{\Dr}\right) \right) \bigg]^{n_K} \notag\\
 & \leq \max_{\hat x_1 \in \Dr} \exp \left( {-P^{\hat x_1}(T_{\hat x} < T_{\DRc}) n_K} \right) \notag \\
 & \leq e^{ -\left[ \left( \frac{\log(R/r) + O(r^{-1/4})}{\log(R)} \right) \left( 1 + O(\log(R)^{-1}) \right) \right] n_K } 
 \leq e^{ -(1 - \de)\frac{\pi_{\Gamma} b (\log K)^2}{2 \log (R/r)} \left( \frac{\log(R/r)}{\log(R)} \right)} \notag\\
 & = e^{ -(1 - \de)\frac{\pi_{\Gamma} b (\log K)^2}{2 \log (R)}} \leq e^{ -(1 - \de)\pi_{\Gamma} b(\log K)/2} \leq K^{-(1 - \de)\pi_{\Gamma} b/2}. \label{eq:BoundSum1}
\end{align}
The second probability in \eqref{eq:SplitLateExcProbs} is bounded above by 
\eqref{eq:lemBoundExcTimesBelow}, 
\begin{align}
P^{\hat x_0}\left(\sum_{j=0}^{n_K} \tau^{(j)} \geq b(K \log K)^2\right)
 & \leq P^{\hat x_0}\left(\sum_{j=0}^{n_K} \tau^{(j)} \geq (1+\de) n_K \frac{2 K^2 \log(R/r)}{\pi_{\Gamma}} \right) \notag\\
 & \leq e^{-C' (1-\de) \pi_{\Gamma} b (\log(K))^2/\log(\log K)}, \label{eq:BoundSum2} 
\end{align}
for some $C' = C'(\delta) > 0$. \eqref{eq:BoundSum1} and \eqref{eq:BoundSum2} combined with \eqref{eq:SplitLateExcProbs} gives us \eqref{eq:LemUpperBd}. $\qed$

The upper bound of (\ref{eq:LatePointsAlpha}) is as follows: For any $\alpha \in (0,1)$ and $\gamma > 0$, we have by Lemma \ref{lem:HittingTimeTailProb}, that for $\gamma/(2\alpha) > \delta > 0$ small enough,
\begin{align} 
 & P\left( \left| \left\{ \hat x \in \Z^2_K: \frac{\mc{T}_K(\hat x)}{(K \log K)^2} \geq \frac{4\alpha}{\pi_{\Gamma}} \right\} \right|
 \geq K^{2(1-\alpha)+\gamma} \right) \notag\\
 \leq & \, K^{-2(1-\alpha)-\gamma} \, 
   \EV\left( \left| \left\{ \hat x \in \Z^2_K: \frac{\mc{T}_K(\hat x)}{(K \log K)^2} \geq \frac{4\alpha}{\pi_{\Gamma}} \right\} \right| \right) \notag\\
 = & \, K^{-2(1-\alpha)-\gamma} \sum_{\hat x \in \Z^2_K} P\left( \frac{\mc{T}_K(\hat x)}{(K \log K)^2} \geq \frac{4\alpha}{\pi_{\Gamma}}\right) \notag\\
 \leq & \, K^{2 \delta \alpha - \gamma} \underset{K \to \infty}{\longrightarrow} 0. \label{eq:LP(3.26)}
\end{align}

\subsection{Lower bound of late point probabilities} \label{sec:LB}

Fixing $0 < \alpha < 1$, we prove in this section the lower bound of \eqref{eq:LatePointsAlpha}: for any $\delta > 0$, $K_n = e^n n^{3n+\gbar}$, and some universal $n_0(\delta) < \infty$, there exists $f_n(\delta) \to 0$ as $n \to \infty$ such that 
\begin{align}
P\left( \bigg|\left\{ \hat x \in \Z^2_{K_n}: \frac{\mc{T}_{K_n}(\hat x)}{(K_n \log K_n)^2} \geq \frac{4\alpha}{\pi_{\Gamma}}\right\}\bigg| \geq K_n^{2(1-\alpha)-\delta} \right) \, \geq \, 1-f_n(\delta). \notag
\end{align}
The sequence $\{K_n\}_{n \geq n_0}$ covers all integers sufficiently to imply 
\begin{align}
\lim_{m \to \infty} & P\left( \bigg|\left\{ \hat x \in \Z^2_m: \frac{\mc{T}_{m}(\hat x)}{(m \log m)^2} \geq \frac{4\alpha}{\pi_{\Gamma}}\right\}\bigg| \geq m^{2(1-\alpha)-\delta} \right) \, = \, 1. \label{eq:LP(4.1)}
\end{align}


Let $a = 2\alpha$ and fix $\rho < \frac{2-a}{2}$. 
We call a pair $(\hat x, \omega)$ \textbf{$n$-successful} \index{nsucc@$n$-successful} if the path $\omega$ does not skip $\hat x$-bands and has the following excursion and visiting counts (where, recall, $v_k = 3ak^2 \log k$): 
\[N_{n,0}^{\hat x} = 0, \,\, |N_{n,k}^{\hat x} - v_k| \leq k, \,\, \ie \,\, N_{n,k}^{\hat x} \stackrel{k}{\sim} v_k, \,\, k = \rho n, \ldots, n-1.\]
Recall that $\mc{R}_n^{\hat x}$ is the time it takes for $v_n$ excursions from $\hp(D(x,r_{n,n-1}))$ out to $\hp(D(x,r_{n,n})^c_{K_n})$ to complete, and note that $\{N_{n,0}^{\hat x} = 0\}$ $=$ $\{\mc{T}_{K_n}(\hat x) > \mc{R}_n^{\hat x}\}$. 
The next lemma relates the notions of $n$-success and first hitting times.  

\begin{lem} \label{lem:LPLemma4.1}
Let $\mc{S}_n = \{\hat x \in \Z^2_{K_n}: \mc{T}_{K_n}(\hat x) > \mc{R}_n^{\hat x}\}$. Then, for some $c > 0$ and all $n \geq n_0$, 
\begin{equation} \label{eq:LPLemma4.1}
P\left( \bigcup_{\hat x \in \mc{S}_n} \left\{ \frac{\mc{T}_{K_n}(\hat x)}{(K_n \log K_n)^2} \leq \frac{2a}{\pi_{\Gamma}} - \frac{2}{\log n} \right\} \right)
 \leq c^{-1} e^{-cn^2/\log n}.
\end{equation}
\end{lem}

\pf Set $r = r_{n,n-1}$, $R = r_{n,n}$, and $\delta = \frac{\pi_{\Gamma}}{2 a \log n}$. Then $\log(R/r) = 3 \log n$, and by \eqref{eq:lemBoundExcTimesAbove} under $N = v_n = 3an^2 \log n$ excursions, we have that, for some $C > 0$, 
all $n \geq n_0$, and any $\hat x, \hat x_0 \in \Z^2_{K_n}$, 
\begin{align*}
P_{\hat x} & := P^{\hat x_0}\left( \mc{T}_{K_n}(\hat x) \leq \left( \frac{2a}{\pi_{\Gamma}} - \frac{2}{\log n} \right)(K_n \log K_n)^2, \mc{T}_{K_n}(\hat x) > \mc{R}_n^{\hat x} \right)\\
  & \leq P^{\hat x_0}\left( \sum_{j=0}^{v_n} \tau^{(j)} \leq \left( \frac{2a}{\pi_{\Gamma}} - \frac{1}{\log n} \right) K_n^2 (3n \log n)^2 \right)\\
  & \leq P^{\hat x_0}\left( \sum_{j=0}^{v_n} \tau^{(j)} \leq (1 - \delta) v_n \frac{2K_n^2 \log(R/r)}{\pi_{\Gamma}} \right) 
  \leq e^{-C \frac{n^2}{\log n}}.
\end{align*}
Sum over $\hat x \in \Z^2_{K_n}$ and select $c < C/2$ so that $c^{-1} e^{-cn_0^2} \geq 1$ to get (\ref{eq:LPLemma4.1}). $\qed$

Let $Y(n,\hat x)$, \index{Ynx@$Y(n,\hat x)$} $\hat x \in \Z^2_{K_n}$, be the indicator random variable for the event 
\[\{\hat x \mbox{ is } n\mbox{-successful}\} =\{\omega: (\hat x,\omega) \mbox{ is } n\mbox{-successful}\}.\]
In view of Lemma \ref{lem:LPLemma4.1}, we have \eqref{eq:LP(4.1)} (and hence \eqref{eq:LatePointsAlpha}) as soon as we show that, for any $\delta > 0$, all $n$ sufficiently large, there exists a sequence $f_n \to 0$ such that 
\begin{equation} \label{eq:LPLemma4.3}
P\left( \sum_{\hat x \in \Z^2_{K_n}} Y(n,\hat x) \geq K_n^{2-a-\delta} \right) \geq 1 - f_n(\delta).
\end{equation}

First, we state \cite[Lemma 6.1]{BRFreq}, a combinatorial result that will aid us in the proof of Lemma \ref{lem:FPLemma5.2}. 
\begin{lem} \label{lem:FPLemma6.1}
For some $C = C(a) < \infty$ and all $k \geq 2$, $|m - v_{k+1}| \leq k+1$, $|l + 1 - v_k| \leq k$, 
\begin{equation} \label{eq:FPLemma6.1(6.14)}
\frac{C^{-1} k^{-3a-1}}{\sqrt{\log k}} \leq \binom{m + l}{l} \left(\frac{1}{2}\right)^{m + l + 1}
 \leq \frac{C k^{-3a-1}}{\sqrt{\log k}}.
\end{equation}
\end{lem}

\begin{lem} \label{lem:FPLemma5.2} 
Fix $\rho < \rho' < \frac{2-a}{2}$. Then there exists $b \geq 10$ and $q_n \geq r_{n,n}^{-a+o(1_n)}$ \st for all $n$ sufficiently large, uniformly in $\gbar \in [b,b+4]$ and $\hat x \in S_{K_n} := \Z^2_{K_n} \setminus \hp(D(0,r_{n,n}))$, \index{qn@$q_n$}
\begin{equation} \label{eq:FPLemma5.2(5.10)}  
P(\hat x \text{ is } n\text{-successful}) = (1 + o(1_n)) q_n.
\end{equation}
\end{lem}

\pf We start by defining a way to examine excursions on a path.  Let $\tau(1)$ be the time of the first visit to $\hp(\bd D(x,r_{n,n-1})_{s_{n-1}^{n \downarrow}})$ (starting at $\hat 0$, so coming from outside $\hat x$'s levels into $\hat x$'s large level $n-1$), and define $\tau(2)$, $\tau(3)$, $\ldots$ to be the successive hitting times of different elements of $A_n := \bigcup_{k=\rho n}^{n} \hp(\bd D(x,r_{n,k})_{s_k})$ until time $\mc{R}_n^{\hat x}$. We can construct a path $\omega$'s ``history'' as follows: let $m = (m_{\rho n}, \ldots, m_{n-1}, m_n)$, where $m_k$ is the number of upcrossing excursions of $\omega$ (candidate values for $N_{n,k}^{\hat x}$) from level $k-1$, \ie $\hp(\bd D(x,r_{n,k-1})_{s_{k-1}})$, out to level $k$, \ie $\hp(\bd D(x,r_{n,k})_{s_k})$ before $\mc{R}_n^{\hat x}$, and set $|\overline{m}| = 2\sum_{k=\rho n}^n m_k - 1$. Let $\Phi : A_n \mapsto \{\rho n - 1, \ldots, n-1, n\}$ label the points of $A_n$ by their annulus: set $\Phi(\hat y) = k$ if $\hat y \in \hp(\bd D(x,r_{n,k})_{s_k})$. Set $h(\omega, j) = \Phi( \omega(\tau(j))$, the label of the annulus hit at time $\tau(j)$, where $\omega \in \Omega_{\hat x,n,n-1,m_n}^{\rho n-1,\ldots,n}$. 
(Note that, since we are referring to upcrossings here, at level $n-1$ we use the thin band $s_{n-1} = n^4$ rather than the thick band $s_{n-1}^{n \downarrow} = \sqrt{r_{n,n-1}}$, which is reserved for the downcrossing $n \downarrow n-1$.)
Since $\omega \in  \Omega_{\hat x,n,n-1,m_n}^{\rho n-1,\ldots,n}$, $h$ satisfies 
\begin{equation} \label{eq:hProperties} 
  h(\omega, 1) = n-1; \,\,\, |h(\omega,j+1) - h(\omega,j)| = 1, j=1,\ldots,|\overline{m}|-1; \,\,\, h(\omega, |\overline{m}|) = n.
\end{equation}
Let $\mc{H}_n(|\overline{m}|)$ be the collection of all such maps 
\[s:\{1,2,\ldots, |\overline{m}|\} \mapsto \{\rho n - 1, \ldots, n-1, n\}\]
satisfying (\ref{eq:hProperties}) for a given $\omega \in  \Omega_{\hat x,n,n-1,m_n}^{\rho n-1,\ldots,n}$. Note that the number of upcrossings from level $k-1$ to $k$ is 
\[ u(k) := |\{ (j,j+1) : (s(j), s(j+1)) = (k-1, k) \}| = m_k.\]
An upcrossing from $k-1$ to $k$ can only occur before the last upcrossing from $k$ to $k+1$. Hence, the number of ways to partition $u(k)$ upcrossings from $k-1$ to $k$ among and before the $u(k+1)$ upcrossings from $k$ to $k+1$ is 
\[\binom{u(k+1) + u(k) - 1}{u(k)},\]
the number of ways to partition $u(k)$ identical objects into $u(k+1)$ sets. Since the mapping $s$ is in one-to-one correspondence with the relative ordering of all its upcrossings, we have 
\[ |\mc{H}_n(\overline{m})| = \prod_{k=\rho n}^{n-1} \binom{m_{k+1} + m_k - 1}{m_k}.\]
Let $h|_k$ be the first $k$ coordinates of the sequence $h$. Applying the strong Markov property at the times $\tau(1)$, $\tau(2)$, ..., $\tau(|\overline{m}|-1)$, we have, uniformly for $s \in \mc{H}_n(\overline{m})$ and $\hat x \in S_{K_n}$, 
\begin{align}
P( h|_{|\overline{m}|} = s;  \Omega_{\hat x,n,n-1,m_n}^{\rho n-1,\ldots,n} & ; \mc{T}_{K_n}(\hat x) > \tau(|\overline{m}|) )  = \prod_{k=\rho n}^{n} a_k^{m_k} b_k^{m_k}, \label{eq:ExcursionsAvoidProb0}
\end{align}
where $a_l$ and $b_l$ are described below.

We wish to examine the probabilities of excursions between annuli. For the outermost level, from level $n$ (\ie the $\hat x$-band of width $s_n = n^4$ at radius $r_{n,n}$), the probability that the toral walk crosses back down to $r_{n,n-1}$ via the thick $\hat x$-band (which is of width $s_{n-1}^{n \downarrow} = \sqrt{r_{n,n-1}}$, unlike all other bands) can be estimated by the bound below \eqref{eq:FP(2.71)}. Uniformly for $\hat w \in \hp(\bd D(x,r_{n,n})_{s_n})$, and for large enough $n$, there exists $c, c' > 0$ \st 
\begin{align}
b_{n} & = P^{\hat w} \left(T_{\hp(D(x,r_{n,n-1}+s_{n-1}^{n \downarrow}))} = T_{\hp(\bd D(x,r_{n,n-1})_{s_{n-1}^{n \downarrow}})} \right) \label{eq:bn0}  \\
 & = 1 - P^{\hat w} \left(T_{\hp(D(x,r_{n,n-1}))} < T_{\hp(\bd D(x,r_{n,n-1})_{s_{n-1}^{n \downarrow}})} \right) \notag\\
 & \geq 1 - c r_{n,n-1} ^{2} \log^2 (r_{n,n-1}) r_{n,n-1}^{-M/2} \notag\\
 & \geq 1 - c r_{n,n-1}^{2-M/2} \log(r_{n,n-1})^2 \notag\\
 & \geq 1 - c' r_{n,n-1}^{-\beta} n^2 (\log n)^2 \notag\\
 &  \geq 1 - c' e^{-\beta n} n^{-3\beta(n-1) + 2} (\log n)^2 = 1 + o(n^{-4}). \notag
\end{align}

From the innermost level $\rho n - 1$, applying \eqref{eq:FP(2.48)}, we will avoid visiting $\hat x$ and cross back up to level $\rho n$ via its $s_{\rho n} = n^4$-band, uniformly in $\hat w \in \hp(\bd D(x,r_{n,\rho n - 1})_{s_{\rho n-1}})$, with probability 
\begin{align}
a_{\rho n} & = P^{\hat w} \left(T_{\hp(D(x,r_{n,\rho n})^c_K)} < T_{\hat x}; \, T_{\hp(D(x,r_{n,\rho n})^c_K)} = T_{\hp(\bd D(x,r_{n,\rho n})_{s_{\rho n}})} \right) \notag \\
 & = 1 - \frac{ \log \left( \frac{ r_{n,\rho n}}{r_{n,\rho n - 1}} \right) + O(r_{n,\rho n - 1}^{-1/4})}{ \log r_{n,\rho n} }\left( 1 + O((\log r_{n,\rho n})^{-1})\right) + o(n^{-8}) \label{eq:arhon} \\
 & = 1 - \frac{3 \log n + o(e^{-n/4})}{n + 3\rho n \log n} (1 + O((\rho n \log n)^{-1})) + o(n^{-8}) \notag \\
 & = 1 - \frac{1}{\rho n} + O\left( (\rho n^2 \log n)^{-1} \right). \notag 
\end{align}
For the middle levels, set $a_l$ to the probability in \eqref{eq:FP6.7analog} for upcrossings for $l=\rho n,\ldots,n$, 
and $b_l$ to \eqref{eq:FP6.8analog} for downcrossings:
\begin{align}
 a_l, \, b_l = \frac{1}{2} + o(n^{-4}), \,\, l=\rho n - 1,\ldots,n-1.  \label{eq:albl}
\end{align}
By \eqref{eq:bn0}, \eqref{eq:arhon}, and \eqref{eq:albl}, \eqref{eq:ExcursionsAvoidProb0} reduces to 
\begin{align}
\prod_{k=\rho n}^{n} a_k^{m_k} b_k^{m_k} & = a_{\rho n}^{m_{\rho n}} b_n^{m_n} \prod_{k=\rho n}^{n-1} a_{k+1}^{m_{k+1}} b_{k}^{m_{k}} \label{eq:ExcursionsAvoidProb} \\
 & = a_{\rho n}^{m_{\rho n}} (1 + o(n^{-4}))^{m_n} \left( \frac{1}{2} + o(n^{-4}) \right)^{|\overline{m}| - m_{\rho n} - m_n + 1} \notag
\end{align}
since $\displaystyle{\sum_{k=\rho n}^{n-1} (m_k + m_{k+1}) = |\overline{m}| - m_{\rho n} - m_n + 1}$. Factoring $\frac{1}{2}$ from the main terms and combining reduces this probability to 
\[ a_{\rho n}^{m_{\rho n}} \left(1 + o(n^{-4})\right)^{|\overline{m}| - m_{\rho n} + 1} \prod_{k=\rho n}^{n-1} \left(\frac{1}{2}\right)^{m_k + m_{k+1}} .\]
Uniformly in $|\overline{m}|$, we have $\left(1 + o(n^{-4})\right)^{|\overline{m}| - m_{\rho n} + 1} = 1 + o(1_n)$.
Finally, for large enough $n$, uniformly in $m_{\rho n} \stackrel{\rho n}{\sim} v_{\rho n}$, and since $a_{\rho n}, \rho \leq 1$, we can bound the term $a_{\rho n}^{m_{\rho n}}$ below: 
\begin{align*}
a_{\rho n}^{m_{\rho n}} & \geq \left(1 - \frac{1}{\rho n} + O((\rho n^2 \log n)^{-1}) \right)^{3a(\rho n)^2 \log(\rho n) + \rho n} \\
 & \geq e^{-3a\rho n \log(\rho n) + O(1)} \geq e^{c} (\rho n)^{3 \rho n(-a)} \\
 & \geq e^{n(-a + o(1_n))} n^{3 \rho n(-a + o(1_n))} \geq r_{n,\rho n}^{-a + o(1_n)}.
\end{align*}
All combined, this yields the exact-history $s$, not-skipping-$\hat x$-bands probability bound 
\begin{align} 
P( h|_{|\overline{m}|} & = s; \Omega_{\hat x,n,n-1,m_n}^{\rho n-1,\ldots,n}; \mc{T}_{K_n}(\hat x) > \tau(|\overline{m}|) ) \notag\\
 & \geq (1 + o(1_n)) r_{n,\rho n}^{-a + o(1_n)} \prod_{k=\rho n}^{n-1}\left(\frac{1}{2}\right)^{m_k + m_{k+1}}. \label{eq:ExcursionsAvoidProb2}
\end{align}

Taking $m_n = v_n = 3an^2 \log n$ and summing over all possible maps $s$ for each possible path $\omega$ gives us 
\begin{align} 
P(\hat x \mbox{ is } n\mbox{-successful}) 
 & = (1 + o(1_n)) \, q_n, \label{eq:qnx1}
\end{align}
which, by \eqref{eq:ExcursionsAvoidProb2}, is \eqref{eq:FPLemma5.2(5.10)} for 
\begin{align} \label{eq:qn2}
q_n \geq r_{n,\rho n}^{-a + o(1_n)} \sum_{\substack{m_{\rho n}, \ldots, m_{n-1}\\|m_k-v_k| \leq k}} \prod_{k=\rho n}^{n-1} \binom{m_{k+1} + m_k - 1}{m_k} \left(\frac{1}{2}\right)^{m_k + m_{k+1}}.
\end{align}
Note that $q_n$ does not depend on $\hat x$. By \eqref{eq:FPLemma6.1(6.14)}, there exists $C, C' < \infty$ independent of $k$ such that, uniformly in $m_k \stackrel{k}{\sim} v_k$ and $m_{k+1} \stackrel{k+1}{\sim} v_{k+1}$,
\begin{align} \label{eq:chooseSqueeze}
\frac{C' k^{-3a-1}}{\sqrt{\log k}} \geq \binom{m_{k+1} + m_k - 1}{m_k} \left(\frac{1}{2}\right)^{m_{k} + m_{k+1}}
 \geq \frac{C k^{-3a-1}}{\sqrt{\log k}}.
\end{align}
Since there are $2l+1$ positive terms for each $l$ \st $m_l \stackrel{l}{\sim} v_l$, the sum in \eqref{eq:qn2} is a sum of $\prod_{l=\rho n}^{n-1} (2l+1)$ terms; each of these terms is a product of $(1-\rho)n$ factors, each of the form $\binom{m_{l+1} + m_l - 1}{m_l} \left(\frac{1}{2}\right)^{m_l + m_{l+1}}$. Thus, using \eqref{eq:chooseSqueeze} and some $C_1, C_1' < \infty$, we can bound the sum in \eqref{eq:qn2} by 
\begin{align} 
\prod_{k=\rho n}^{n-1} \frac{C_1' k^{-3a}}{\sqrt{\log k}} & \geq\prod_{l=\rho n}^{n-1} (2l + 1) \prod_{k=\rho n}^{n-1} \frac{C' k^{-3a-1}}{\sqrt{\log k}} \notag\\
& \geq \sum_{\substack{m_{\rho n}, \ldots, m_{n-1}\\|m_l-v_l| \leq l}} \prod_{k=\rho n}^{n-1} \frac{C k^{-3a-1}}{\sqrt{\log k}} \geq \prod_{l=\rho n}^{n-1} (2l + 1) \prod_{k=\rho n}^{n-1} \frac{C k^{-3a-1}}{\sqrt{\log k}} \label{eq:chooseSqueezeSum} \\
 & \geq \prod_{k=\rho n}^{n-1} \frac{C_1 k^{-3a}}{\sqrt{\log k}} \geq (1-\rho) n C_1^{(1-\rho)n} n^{3(1-\rho)n(-a)} \left(\prod_{k=\rho n}^{n-1} \log k\right)^{-1/2}. \notag
\end{align}
It is obvious that a constant $c$ is $n^{o(1_n)}$, and $n^c$ is $(n^n)^{o(1_n)}$ 
for any fixed $c>0$. Hence,
\begin{equation} \label{eq:killSmallThings}
(1-\rho) n C_1^{(1-\rho)n} = (n^n)^{o(1_n)} = r_{n,n}^{o(1_n)}.
\end{equation}
Next, $n^{3(1-\rho)n(-a)}$ combined with $r_{n,\rho n}^{-a + o(1_n)}$ yields 
\begin{equation} \label{eq:rnnComesAlive}
r_{n,\rho n}^{-a + o(1_n)} n^{3(1-\rho)n(-a)} = (e^n n^{3 \rho n})^{-a+o(1_n)} (n^{3(1-\rho)n})^{-a} = r_{n,n}^{-a + o(1_n)}.
\end{equation}
Finally, 
\begin{align} 
\left(\prod_{k=\rho n}^{n-1} \log k\right) = n^{nx} 
 \implies & x = \frac{\log \left(\prod_{k=\rho n}^{n-1} \log k\right)}{n \log n} \leq \frac{(1-\rho)n\log \log n}{n \log n} \to 0 \notag\\
 \implies & \left(\prod_{k=\rho n}^{n-1} \log k\right)^{-1/2} = r_{n,n}^{o(1_n)}. \label{eq:biglogbound}
\end{align}
Merging \eqref{eq:qn2}-\eqref{eq:biglogbound} results in $q_n \geq r_{n,n}^{-a + o(1_n)}$. $\qed$

For a given $n$, define \index{lxy@$l(\hat x,\hat y)$}
\[l(\hat x,\hat y) := \max\{m \in \{0,1,2,\ldots, n\}: \hp(D(x,r_{n,m})) \cap \hp(D(y,r_{n,m})) = \emptyset\}\]
to be the largest radius index (up to $n$) of discs centered at $\hat x$ and $\hat y$ that do not intersect. We now show that the covariance of $Y(n,\hat x)$ between pairs of points depends on how far apart they are, based on this measurement.

\begin{lem} \label{lem:FPLemma5.2(5.12)}  
Fix $\varepsilon > 0$. Then there exists $b \geq 10$ and $C = C(b, \varepsilon)< \infty$ 
 \st for all $n$ and $\hat x, \hat y \in S_{K_n}$, 
\begin{align} \label{eq:FPLemma5.2(5.12)} 
\EV( Y(n,\hat x) Y(n,\hat y) ) \leq \left\{\begin{array}{ll}
C^{n-l} q_n^2 n^b \left( \frac{r_{n,n}}{r_{n,l}} \right)^{a+\varepsilon} & \,\, \rho' n \leq l(\hat x,\hat y) < n, \\
(1 + o(1_n)) q_n^2
 & \,\, l(\hat x,\hat y) = n.
\end{array}\right.
\end{align}
\end{lem}

\pf First, note that, using the index set $M_l := \{l, l+1, \ldots, n-1\}$, \index{Ml@$M_l$} the same analysis at the end of the proof of Lemma \ref{lem:FPLemma5.2} yields, for any $l \geq \rho n$, uniformly in $\hat x \in S_{K_n}$, $\gbar$, and $m_k \leq 3k^2 \log k + k$, 
\begin{equation} \label{eq:LP(5.9)}
P(N_{n,k}^{\hat x} = m_k, k \in M_l) = (1 + o(1_n)) \prod_{k=l}^{n-1} \binom{m_{k+1} + m_k - 1}{m_k} \left( \frac{1}{2} \right)^{m_k + m_{k+1}}.
\end{equation}
Recall that $v_k = v_k(a) = 3ak^2 \log k$ and $N \stackrel{k}{\sim} v_k$ if $|N - v_k| \leq k$ for $\rho n \leq k < n$ and $N=0$ if $k=0$. We first note that, for $\rho' n \leq l(\hat x,\hat y)<n$, 
$2r_{n,l+1} + 2 \geq d(\hat x, \hat y) \geq 2r_{n,l} + 2$. Thus, there are, for some constants $C_{n,k} \approx 4\pi$, 
\begin{equation} \label{eq:DecomposeZ2KnByDistance}
|\{y: l(\hat x,\hat y) = l\}| = C_{n,l+1}(r_{n,l+1}^2 - r_{n,l}^2). 
\end{equation}
Since $r_{n,l+2} - r_{n,l} \gg r_{n,l+1}$, it is easy to see that 
\[l = l(\hat x,\hat y)<n \implies \hp(D(y,r_{n,l}')) \cap \hp(\bd D(x,r_{n,k})_{s_k}) = \emptyset\]
for $k \neq l+1$ (the thick band at $k=n-1$ also satisfies this). Replacing hereafter $l$ with $l \land n-3$, it follows that for $k \neq l+1,l+2$, the events 
 $\{N_{n,k}^{\hat x} \stackrel{k}{\sim} v_k\}$ are measurable with respect to the $\sigma$-algebra $\mc{G}_{n,l \downarrow l-1}^{\hat y}$ (defined before Lemma \ref{lem:FPLemma8.2}), since the excursions outside these bands depend (up to error term) only on their beginning and end points. 

\begin{figure}[!ht]
  \centering
    \includegraphics[width=6in]{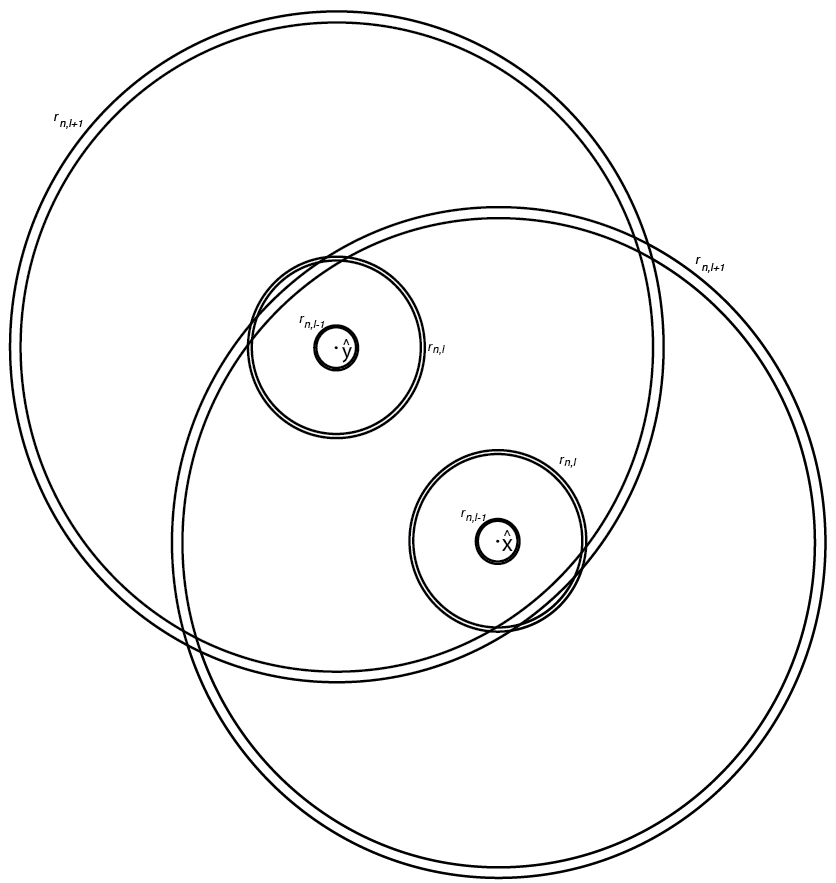}
  \parbox{4in}{
  \caption{An example of $l(\hat x,\hat y) = l$ where levels $l$ and $l+1$ have nonempty intersection}
  \label{fig:level_l_vs_level_lplus1}}
\end{figure}

Slightly rewriting the notation of Lemma \ref{lem:FPCor8.3}, define the set of $\hat y$-\emph{faithful} \index{faithful} paths for the set of indices $A$, \index{GnlyA@$\Gamma_n^{\hat y}(A)$}
\[\Gamma_{n}^{\hat y}(A) := \{ N_{n,i}^{\hat y} \stackrel{i}{\sim} v_i; \,\, i \in A\} \cap \Omega_{y,n,l,v_l+l}^{A},\]
to be the set of paths with $n$-successful $\hat y$-excursion counts on the levels of the indices of $A$. Using the index set $J_l = \{0, \rho n,\ldots,l-1\}$, \index{Jl@$J_l$} we collect all the pertinent inner-level $\hat y$-based excursions, and with the index set $I_l = \{0, \rho n, \ldots, l, l+3,\ldots,n-1\}$,  we combine the inner- and outer-level $\hat x$-faithful excursion paths, skipping the two levels where $\hat x$ and $\hat y$'s annuli cross (causing a jump in their $n$-success covariance).

Note that $\Gamma_n^{\hat x}(I_l) \in \mc{G}_{n,l \downarrow l-1}^{\hat y}$ (it skips the two levels in question). Then we have that  
\[\{\hat x \mbox{ and } \hat y \mbox{ are } n \mbox{-successful}\} \subset \Gamma_{n}^{\hat x}(I_l) \cap \Gamma_{n}^{\hat y}(J_{l+1}). \]



Recall that, if $B \in \mc{G}$, $P(A \cap B | \mc{G}) = P(A | \mc{G}) 1_B$. Applying \eqref{eq:FP(8.8)}, and focusing on level $l$, for some universal constant $C_3 < \infty$, 
\begin{align} 
P(\hat x \text{ and } \hat y \text{ are } n \text{-successful})
 & \leq \sum_{m_l \stackrel{l}{\sim} v_l} \EV \left( P \left( \Gamma_{n}^{\hat y}(J_l) | N_{n,l}^{\hat y} = m_l, \mc{G}_{n,l \downarrow l-1}^{\hat y}\right); \, \Gamma_n^{\hat x}(I_l)  \right) \notag \\
 & \leq C_3 P(\Gamma_n^{\hat x}(I_l)) \sum_{m_l \stackrel{l}{\sim} v_l} P\left( \Gamma_{n}^{\hat y}(J_l) | N_{n,l}^{\hat y} = m_l \right). \label{eq:LP(5.10)}
\end{align}

Using Lemma \ref{lem:FPCor8.3}, for some universal constant $0 < C_4 < \infty$, 
\begin{align} 
(1 + o(1_n)) q_n & = P(y \text{ is } n\text{-successful}) \label{eq:LP(5.11)}\\
 & = \sum_{m_l \stackrel{l}{\sim} v_l} \EV \left( P\left( \Gamma_{n}^{\hat y}(J_l) | N_{n,l}^{\hat y} = m_l, \mc{G}_{n,l \downarrow l-1}^{\hat y}  \right); \, 
 N_{n,l}^{\hat y} = m_l, \Gamma_n^{\hat y}(M_{l+1}) \right) \notag\\
 & \geq C_4 \sum_{m_l \stackrel{l}{\sim} v_l} P\left( N_{n,l}^{\hat y} = m_l, \Gamma_n^{\hat y}(M_{l+1}) \right) \times P(\Gamma_n^{\hat y}(J_l) | N_{n,l}^{\hat y} = m_l). \notag
\end{align}



Hence, by \eqref{eq:LP(5.9)} and \eqref{eq:chooseSqueeze}, 
for some universal $C_5 < \infty$, 
\begin{equation} \label{eq:LP(5.12)} 
\sum_{m_l \stackrel{l}{\sim} v_l} P\left( \Gamma_{n}^{\hat y}(J_l) | N_{n,l}^{\hat y} = m_l \right) \leq C_5^{n-l} q_n l \left( \prod_{k=l}^{n-1} k^{3a} \sqrt{\log k} \right).
\end{equation}
Similarly, using Lemma \ref{lem:FPCor8.3}, 
\begin{align} 
P(\Gamma_n^{\hat x}(I_l)) & \leq \sum_{m_l \stackrel{l}{\sim} v_l} \EV \left( P\left( \Gamma_n^{\hat x}(J_l) | N_{n,l}^x = m_l, \mc{G}_{n,l \downarrow l-1}^x \right); \, \Gamma_n^{\hat x}(M_{l+3}) \right) \label{eq:LP(5.13)} \\
 & \leq C_6 P(\Gamma_n^{\hat x}(M_{l+3}) \sum_{m_l \stackrel{l}{\sim} v_l}  P\left( \Gamma_n^{\hat x}(J_l) | N_{n,l}^x = m_l \right). \notag
\end{align}
Comparing \eqref{eq:LP(5.13)} and \eqref{eq:LP(5.11)}, and applying \eqref{eq:LP(5.9)} and \eqref{eq:chooseSqueeze}
again, we get 
\begin{align} \label{eq:LP(5.14)}
P(\Gamma_n^{\hat x}(I_l)) & \leq C_7 l \left( \prod_{k=l}^{l+2} k^{3a} \sqrt{\log k} \right) q_n.
\end{align}

Combining \eqref{eq:LP(5.10)}, \eqref{eq:LP(5.12)},  and \eqref{eq:LP(5.14)} proves \eqref{eq:FPLemma5.2(5.12)} for $l(\hat x,\hat y) < n$.

Finally, we deal with those pairs far apart. For most pairs ($K_n^2(K_n^2 - C_{n,n}r_{n,n}^2)$ pairs for some $C_{n,n} \approx 4\pi$ of them), we have $l(\hat x,\hat y) = n$. For these, the event $\{\hat x$ is $n$-successful$\}$ is  $\mathcal{G}_{n,n \downarrow n-1}^{\hat y}$-measurable, so by Lemma \ref{lem:FPCor8.3}, 
\begin{align} 
\EV( Y(n,\hat x) Y(n,\hat y)) & = P(\hat x \mbox{ and } \hat y \mbox{ are } n\mbox{-successful}) \notag\\
 & = \EV( P( \hat y \mbox{ is } n\mbox{-successful} \, |\, \mathcal{G}_{n,n \downarrow n-1}^y); \hat x \mbox{ is } n\mbox{-successful}) \label{eq:EndThisMadness}\\
 & \leq (1 + O(n^{-1} (\log n)^2))(1 + o(1_n)) q_n^2 = (1 + o(1_n)) q_n^2. \qed \notag
\end{align}

We can now prove Theorem \ref{thm:LatePoints}. 

Let \index{vl@$V_l$}
\[ V_l := \sum_{x,y \in S_{K_n}, l(\hat x, \hat y) = l} \EV(Y(n,\hat x), Y(n,\hat y)), \,\,\, l=0,1,\ldots,n.\]
Since, by \eqref{eq:FPLemma5.2(5.10)}, considering the sum $\displaystyle{W_n := \sum_{\hat x \in S_{K_n}} Y(n,\hat x)}$, the number of $n$-successful points $\hat x$, 
\[ \EV(W_n) = \EV \left( \sum_{x \in S_{K_n}} Y(n,\hat x) \right) = (1 + o(1_n)) K_n^2 q_n \geq K_n^{2-a+o(1_n)},\]
recall the Paley-Zygmund inequality (\cite[Lemma 14.8.2]{MR}): since $W_n \in L^2(\Omega)$, for any $0 < \lambda_n < 1$, we have 
\begin{align} \label{eq:PZ}
P(W_n \geq \lambda_n \EV(W_n)) \geq (1 - \lambda_n)^2 \frac{\EV(W_n)^2}{\EV(W_n^2)}.
\end{align}
By \eqref{eq:PZ}, \eqref{eq:LPLemma4.3} will follow from the bottom half of (7.34) and 
\begin{align} \label{eq:LP(4.7)}
\EV(W_n^2) = \sum_{l=0}^{n-1} V_l \leq o(1_n) K_n^4 q_n^2.
\end{align}
To obtain this bound, first note that the definition of $l(x,y)$ implies that $ d(\hat x, \hat y)<2r_{n,l(x,y)+1} + 2 $. Hence, on $\Z^2_{K_n}$ there are at most $C_{0} r_{n,l+1}^2$ points $\hat y \in \hp(D(x, r_{n,l+1}))$ (from here on, $C_m$ are constants independent of $n$). 
Since $2 \rho' < 2 - a$, there exists $C_1 < \infty$ such that the covariances on the inner levels sum to 
\begin{align} 
\sum_{l=0}^{\rho' n-1} V_l & \leq \sum_{\hat x, \hat y \in \Z^2_{K_n},d(\hat x, \hat y)\leq  2 r_{n,\rho' n}} \EV(Y(n,\hat x) Y(n,\hat y)) \label{LP(4.8)} \\
 & \leq \sum_{\hat x, \hat y \in \Z^2_{K_n},d(\hat x, \hat y)\leq  2 r_{n,\rho' n}} \EV(Y(n,\hat x)) \leq C_1 q_n K_n^2 r_{n,\rho' n}^2 \leq o(1_n) K_n^4 q_n^2. \notag
\end{align}
Choose $\varepsilon > 0$ such that $2 - a - \varepsilon > 0$ and fix $l \in [\rho' n, n)$. Then, by \eqref{eq:FPLemma5.2(5.12)}, the outer-level covariances are bounded by 
\begin{equation} \label{eq:LP(4.8a)}
V_l \leq C_2 K_n^2 r_{n,l+1}^2 q_n^2 n^b C^{n-l} \left( \frac{r_{n,n}}{r_{n,l}} \right)^{a+\varepsilon},
\end{equation}
which leads to the overall upper-level covariance bound 
\begin{align} 
\sum_{l=\rho' n}^{n-1} V_l & \leq C_2 K_n^2  q_n^2 n^b \sum_{l=\rho' n}^{n-1} C^{n-l} r_{n,l+1}^2 \left( \frac{r_{n,n}}{r_{n,l}}  \right)^{a+\varepsilon} \notag \\
 & = C_2 K_n^4 q_n^2 n^{-2\gbar + b + 6} \sum_{l=\rho' n}^{n-1} C^{n-l} \left( \frac{r_{n,l}}{r_{n,n}}  \right)^{2-a-\varepsilon}  \label{eq:LP(4.9)} \\
 & \leq C_2 K_n^4 q_n^2 n^{-2} \sum_{j=1}^{n} C^{j} r_{n,j}^{-(2-a-\varepsilon)}. \notag 
\end{align}
Combining \eqref{LP(4.8)} and \eqref{eq:LP(4.9)} we get \eqref{eq:LP(4.7)}, which proves \eqref{eq:LPLemma4.3} and thus \eqref{eq:LP(4.1)}. $\qed$

Finally, we prove the cover time result, Corollary \ref{thm:CoverTime}.

\pf The lower bound \eqref{eq:LP(4.1)} implies that, for any $\alpha \in (0,1)$, $\alpha$-late points exist with positive probability. 
As $\alpha \uparrow 1$, we have that $\frac{\mc{T}_{cov}(\Z^2_K)}{(K \log K)^2} \geq \frac{4}{\pi_{\Gamma}}$ in probability as $K \to \infty$. 

For the upper bound, we modify the argument of \eqref{eq:LP(3.26)} to approach from above, \ie $\alpha \downarrow 1$, to show that, as $K \to \infty$, we have no late points beyond $\alpha = 1$ after the expected cover time $\frac{4}{\pi_{\Gamma}}(K \log K)^2$. Define, for any $\alpha > 0$, the cover time event 
\[ A_{\alpha}^K := \left\{ |\mc{L}_K(\alpha)| = \left| \left\{ \hat x \in \Z^2_K: \frac{\mc{T}_K(\hat x)}{(K \log K)^2} \geq \frac{4\alpha}{\pi_{\Gamma}} \right\} \right| \geq K^{0} = 1 \right\} = \left\{\frac{\mc{T}_{cov} (\Z^2_K)}{(K \log K)^2} > \frac{4\alpha}{\pi_{\Gamma}} \right\}. \]
For any $\delta > 0$, set 
$b = \frac{4\alpha}{\pi_{\Gamma}} = \frac{4+\delta}{\pi_{\Gamma}(1-\delta)}$ (so that $\alpha = \frac{4+\delta}{4(1-\delta)} > 1$); Lemma \ref{lem:HittingTimeTailProb} and \eqref{eq:LP(3.26)} yield 
\begin{align*} 
P(A^K_{\alpha}) = P(|\mc{L}_K(\alpha)| \geq 1)
 & = P\left( \left| \left\{ \hat x \in \Z^2_K: \frac{\mc{T}_K(\hat x)}{(K \log K)^2} \geq b \right\} \right| \geq 1 \right) \notag\\
 & \leq \EV\left( \left| \left\{ \hat x \in \Z^2_K: \frac{\mc{T}_K(\hat x)}{(K \log K)^2} \geq b \right\} \right| \right) \notag\\
 & = \sum_{\hat x \in \Z^2_K} P\left( \frac{\mc{T}_K(\hat x)}{(K \log K)^2} \geq b\right) \notag\\
 & \leq \, K^{2 -(1-\delta)\pi_{\Gamma} b/2} = K^{-\delta/2} \underset{K \to \infty}{\longrightarrow} 0. \qed
\end{align*}


\section{Open Problems}

We have given the asymptotic timing of a large class of infinite-range symmetric random walks on the two-dimensional torus. Some open problems to extend this work are: 

\begin{itemize}
\item Analyze the neighborhoods and pairs of late points mentioned in \cite[Theorems 1.2 and 1.3]{DPRZ2006}. How is the spacing of $\alpha$-late point pairs on $\Z^2_K$ affected by jumping walks? 
\item Examine the structure of the frequent points on the lattice torus. 
\item \cite{DPRZ2006} suggests that its nearest-neighbor results may be extended to the planar Weiner sausage on the two-dimensional torus $\mathbb{T}^2$. We suggest, then, that using this class of jumping walks, this work may be extended to a larger class of ``compound Poisson Weiner sausage links'' on $\mathbb{T}^2$ (for example, a two-dimensional Brownian motion with exponentially-timed jumps). 
\item Check the ratio of late points of $\Z_{K_1} \times \Z_{K_2}$ when limiting the coordinates at different rates and when limiting to the infinite cylinder $\Z^2 \times \Z_K$ for fixed $K$.
\item Find tight bounds for $\hGDn(\hat x,\hat x)$, the external toral Green's function, along with annulus Green's functions on the plane and torus and expected hitting times of these discs and annuli, and prove a full exterior toral Harnack inequality.
\item Give computational rates of convergence for the number of late points, given $\alpha$ and $p_1$.
\end{itemize}


\singlespacing
\cleardoublepage

\end{document}